\documentclass[12pt,twoside,a4paper]{article}
\usepackage{amsmath,amssymb,latexsym,theorem,bbm,natbib,epsfig,color}
\usepackage{multirow,graphicx,array,rotating,epstopdf,easybmat}  
\newfont{\msa}{msam10 scaled\magstep1}
\newfont{\ssmsa}{msam9}

\newcommand{\tr}{\mathop{\hbox{\rm tr}}}

\numberwithin{equation}{section}

\newcommand{\proofend}{\hfill$\square$}

\newtheorem{thm}{Theorem}[section]

\newtheorem{lem}[thm]{Lemma}

\theorembodyfont{\rm}

\newtheorem{ex}[thm]{Example}

\title{D-optimal designs for complex Ornstein-Uhlenbeck processes}

\author{{\sc S\'andor Baran$^{1}$}, {\sc 
    Csilla Sz\'ak-Kocsis$^{1}$} and {\sc Milan Stehl\'\i k$^{2,3}$}\\ 
         $^1$Faculty of Informatics, University of Debrecen\\
         Kassai \'ut 26, H-4028 Debrecen, Hungary \\
         $^2$ Department of Applied Statistics and  Linz Institute of 
         Technology, JKU\\
         Altenberger Stra{\ss}e 69, A-4040, Austria\\
         $^3$ Institute of Statistics, University of Valpara\'\i so\\
         Gran Breta\~na 1111, Valpara\'\i so, Chile}

\date{}

\begin{document}
\pagestyle{myheadings}

\maketitle

\begin{abstract}
Complex Ornstein-Uhlenbeck (OU) processes have various applications in statistical modelling. They play role e.g. in the description of the motion of a charged test particle in a constant magnetic field or in the study of rotating waves in time-dependent reaction diffusion systems, whereas Kolmogorov used such a process to model the so-called Chandler wobble, small deviation in the Earth's axis of rotation.  In these applications parameter estimation and model fitting is based on discrete observations of the underlying stochastic process, however, the accuracy of the results strongly depend on the observation points.

This paper studies the properties of D-optimal designs for estimating the parameters of a complex OU process with a trend.  We show that in contrast with the case of the classical real OU process, a D-optimal design exists not only for the trend parameter, but also for joint estimation of the covariance parameters, moreover, these optimal designs are equidistant.

\end{abstract}

\section{Introduction}
  \label{sec:sec1}
Random processes have various applications in statistical modelling in different areas of science such as physics, chemistry, biology or finance,   where one usually cannot observe continuous trajectories. In these situations parameter estimation and model fitting is based on discrete observations of the underlying stochastic process, however, the accuracy of the results strongly depend on the observation points. The theory of optimal experimental designs, dating back to the late 50s of the twentieth century \citep[see e.g.][]{hoel58, kiefer59}, deals with finding design sets \ $\boldsymbol\xi = \{t_1, t_2 , \ldots , t_n \}$ \ of distinct time points (or locations in space) where the process under study is observed, which are optimal according to some previously specified criterion \citep{muller}. In parameter estimation problems the most popular criteria are based on the Fisher information matrix (FIM) of the observations. D-, E- and T-optimal designs maximize the determinant, the smallest eigenvalue and the trace of the FIM, respectively, an A-optimal design minimizes the trace of the inverse of the FIM  \citep[for an overview see][]{puk}, whereas K-optimality refers to the minimization of the condition number of the FIM \citep[see e.g.][]{yz,baran}. In the last decades information based criteria have intensively been studied and despite the well developed theory for uncorrelated setup \citep[see e.g.][]{silvey} only recently the more difficult correlated situation has been addressed \citep{dpz,dpz2}. 

In the present paper we derive D-optimal designs for parameter estimation of complex (or vector) Ornstein-Uhlenbeck (OU) processes with trend \citep[see e.g.][]{arato}, defined in detail in Section \ref{sec:sec2}.
A complex OU process describes e.g. the motion of a charged test particle in a constant magnetic field \citep{balescu}, it is used in the description of the rotation of a planar polymer \citep{vyh} or in the study of rotating waves in time-dependent reaction diffusion systems \citep{bl,otten}, and it also has several applications in financial mathematics \citep[see e.g.][]{bns}. Further, Kolmogorov proposed to model the so-called Chandler wobble, small deviation in the Earth's axis of rotation \citep{lambeck}, by the model
\begin{equation}
  \label{chandler}
Z(t)=Z_1(t)+iZ_2(t)=m{\mathrm e}^{i2\pi t}+Y(t), \qquad t>0,
\end{equation}
where \ $Z_1(t)$ \ and \ $Z_2(t)$ \ are the coordinates of the deviation of the instantaneous pole from the North Pole and \ $Y(t)$ \ is a complex OU process \citep{aks}.

We remark that the properties of D-optimal designs for classical one-dimensional OU processes have already investigated by \citet{ks} and later by \citet{zba09}, where the authors proved that there is no D-optimal design for estimating the covariance parameter, whereas the D-optimal design for trend estimation is equidistant and larger distances resulting in more information. Later these results were generalized for OU sheets under various sampling schemes \citep{bs15,bss13,bss15}.
\vfill

\section{Complex Ornstein-Uhlenbeck process with a trend}
  \label{sec:sec2}

Consider the complex stochastic process \ $Z(t)=Z_1(t)+iZ_2(t)$, \ defined as
\begin{equation}
  \label{c_mod}
Z(t)=mf(t)+Y(t), \qquad t\geq 0,
\end{equation}
with design points taken from the non-negative half-line \ ${\mathbb R_+}$, \ where \ $m=m_1+im_2, \ m_1,m_2\in{\mathbb R}$, \ $f(t)=f_1(t)+if_2(t)$ \ with \ $f_1(t),f_2(t): {\mathbb R_+}\mapsto {\mathbb R}$ \ and \ $Y(t)=Y_1(t)+iY_2(t), \ t\geq 0$, \ is a complex valued stationary OU process.  The process \ $Y(t)$ \ can be defined by the stochastic differential equation (SDE)
\begin{equation}
  \label{c_sde}
{\mathrm d}Y(t)=-\gamma Y(t){\mathrm d}t+ \sigma {\mathrm d}{\mathcal W}(t), \qquad Y(0)=\xi,
\end{equation}
where \ $\gamma=\lambda -i\omega$ \ with \ $\lambda>0, \ \omega\in {\mathbb R}$, \ $\sigma>0$,  \ ${\mathcal W}(t)={\mathcal W}_1(t)+i{\mathcal W}_2(t), \ t\geq 0$, \ is a standard complex Brownian motion, that is \ ${\mathcal W}_1(t)$ \ and \ ${\mathcal W}_2(t)$ \ are independent standard Brownian motions, and \ $\xi=\xi_1+i\xi_2$, \ where \ $\xi_1$ \ and \ $\xi_2$ \ are centered normal random variables that are chosen according to stationarity \citep{arato}.

Instead of the complex process \ $Y(t)$ \ defined by \eqref{c_sde} one can consider the two-dimensio\-nal real valued stationary OU process \ $\big( Y_1(t),Y_2(t)\big)^{\top}$ \ defined by the SDE
\begin{equation}
  \label{2d_sde}
  \begin{bmatrix} {\mathrm d}Y_1(t) \\ {\mathrm d}Y_2(t)\end{bmatrix} =
  A
   \begin{bmatrix} Y_1(t) \\ Y_2(t)\end{bmatrix}{\mathrm d}t
   +\sigma  \begin{bmatrix} {\mathrm d}{\mathcal W}_1(t) \\ 
     {\mathrm d}{\mathcal W}_1(t)\end{bmatrix}, \quad \text{where} \quad 
   A:=\begin{bmatrix} -\lambda & -\omega \\ \omega & -\lambda \end{bmatrix}
\end{equation}
We remark that in physics \eqref{2d_sde} is called A-Langevin equation, see e.g. \citep{balescu}.
If \ $\big( Y_1(t),Y_2(t)\big)^{\top}$ \ satisfies \eqref{2d_sde} then \ $Y_1(t)+iY_2(t)$ \ is a complex OU process which solves \eqref{c_sde}, and conversely, the real and imaginary parts of a complex OU process form a two-dimensional real OU process satisfying \eqref{2d_sde}. Obviously, \ ${\mathsf E}  Y_1(t)={\mathsf E}  Y_2(t)=0$, \ whereas the covariance matrix function of the process  \ $\big( Y_1(t),Y_2(t)\big)^{\top}$ \ is given by
\begin{equation}
 \label{ou_cov}
{\mathcal R}(\tau):= {\mathsf E}  \begin{bmatrix} Y_1(t+\tau ) \\ Y_2(t+\tau )\end{bmatrix}  \begin{bmatrix} Y_1(t) \\ Y_2(t)\end{bmatrix}^{\top} \!\!= \frac {\sigma^2}{2\lambda} {\mathrm e}^{A\tau}=\frac {\sigma^2}{2\lambda} {\mathrm e}^{-\lambda \tau}\begin{bmatrix} \cos(\omega\tau)  & -\sin(\omega\tau) \\ \sin (\omega\tau) & \cos (\omega\tau) \end{bmatrix}, \quad \tau \geq 0.
\end{equation}
This results in a complex covariance function of the complex OU process \ $Y(t)$ \ defined by \eqref{c_sde} of the form
\begin{equation*}
{\mathcal C}(\tau):={\mathsf E}Y(t+\tau)\overline{Y(t)}=\frac {\sigma^2}{\lambda} {\mathrm e}^{-\lambda \tau}\big( \cos(\omega\tau)+i\sin(\omega\tau)\big ), \qquad \tau\geq 0,
\end{equation*}
behaving like a damped oscillation with frequency \ $\omega$.

In the present study the damping parameter \ $\lambda$, \ frequency \ $\omega$ and standard deviation \ $\sigma$ \  are assumed to be known. However, a valuable direction for future research will be the investigation of models
where these parameters should also be estimated. Note that the estimation of \ $\sigma$ \ can easily be done on the basis of a single realization of the complex process, see e.g. \citet[][Chapter 4]{arato}. Now, without loss of generality, one can set the variances of \ $Y_1(t)$ \ and \ $Y_2(t)$ \ to be equal to one, that is \ $\sigma^2/(2\lambda)=1$, \  which reduces \ ${\mathcal R}(\tau)$ \ to a correlation matrix function. Further results on the maximum-likelihood estimation of the covariance parameters can be found e.g. in \citet{abi}.

\section{D-optimal designs}
   \label{sec:sec3}

Suppose the complex process \ $Z(t)$ \   
is observed in design points \ $0\leq t_1<t_2 < \ldots <t_n$ \ resulting in the $2n$-dimensional real vector \ $\big( Z_1(t_1),Z_2(t_1),Z_1(t_2),Z_2(t_2), \ldots ,Z_1(t_n),Z_2(t_n)\big)^{\top}$\!\!, where
\begin{equation}
  \label{2d_mod}
Z_1(t)=m_1f_1(t)-m_2f_2(t)+Y_1(t), \qquad Z_2(t)=m_2f_1(t)+m_1f_2(t)+Y_2(t).
\end{equation} 

As it has mentioned in the Introduction, a D-optimal design maximizes the determinant of the FIM on the unknown parameters corresponding to the observations. Here we consider optimal designs for estimating the trend parameter \ $m$ \ and the damping parameter \ $\lambda$ \ and frequency \ $\omega$, \ as well.

\subsection{Estimation of the trend parameter}
  \label{subs:subs3.1}

According to the results of e.g. \citet{xia} or \citet{pazman} the FIM on parameter vector \ $(m_1,m_2)^{\top}$ \ based on observations \ $\big\{ \big(Z_1(t_j),Z_2(t_j)\big), \ j=1,2,\ldots ,n\big\}$ \ equals
\begin{equation*}
\mathcal I_{m_1,m_2}(n)=H(n)C(n)^{-1}H(n)^{\top}, 
\end{equation*}
where
\begin{equation*}
H(n):=\begin{bmatrix} f_1(t_1) & f_2(t_1) & f_1(t_2) & f_2(t_2) &\cdots &  f_1(t_n) & f_2(t_n) \\ -f_2(t_1) & f_1(t_1) & -f_2(t_2) & f_1(t_2) &\cdots &  -f_2(t_n) & f_1(t_n) \end{bmatrix},
\end{equation*}
and  \ $C(n)$ \ is the  \ $2n \times 2n$ \ covariance matrix of the observations.

\begin{lem}
  \label{lem1}
The FIM on trend parameters  \ $(m_1,m_2)^{\top}$ \ of the two-dimensional analogue of the complex model \eqref{c_mod} based on the real observation vector \ $\big\{\! \big(Z_1(t_j),Z_2(t_j)\big), \ j\!=\!1,2,\ldots ,n\big\}$ \ equals \ $\mathcal I_{m_1,m_2}(n)= Q(n)\,{\mathbb I}_2$, \ where  \  ${\mathbb I}_k, \ k\in{\mathbb N}$, \ denotes the $k$-dimensional unit matrix and
\begin{align}
  \label{Qn_general}
&Q(n):=f_1^2(t_n)\!+\!f_2^2(t_n)\!+\!\sum_{j=1}^{n-1} \frac 1{1-{\mathrm e}^{-2\lambda d_j}}\bigg(f_1^2(t_j)+f_2^2(t_j)+{\mathrm e}^{-2\lambda d_j} \big(f_1^2(t_{j+1})+f_2^2(t_{j+1})\big) \\
&+\!2{\mathrm e}^{-\lambda d_j}\!\Big(\!\big(f_1(t_j)f_2(t_{j+1})\!-\!f_2(t_j)f_1(t_{j+1})\big)\!\sin (\omega d_j) \!-\!\big(f_1(t_j)f_1(t_{j+1})\!+\!f_2(t_j)f_2(t_{j+1})\big)\!\cos (\omega d_j)\!\Big)\!\bigg)\!, \nonumber
\end{align}
with \ $d_j:=t_{j+1}-t_j, \ j=1,2,\ldots, n-1$.
\end{lem}

\begin{figure}[t]
\begin{center}
\leavevmode
\hbox{
\epsfig{file=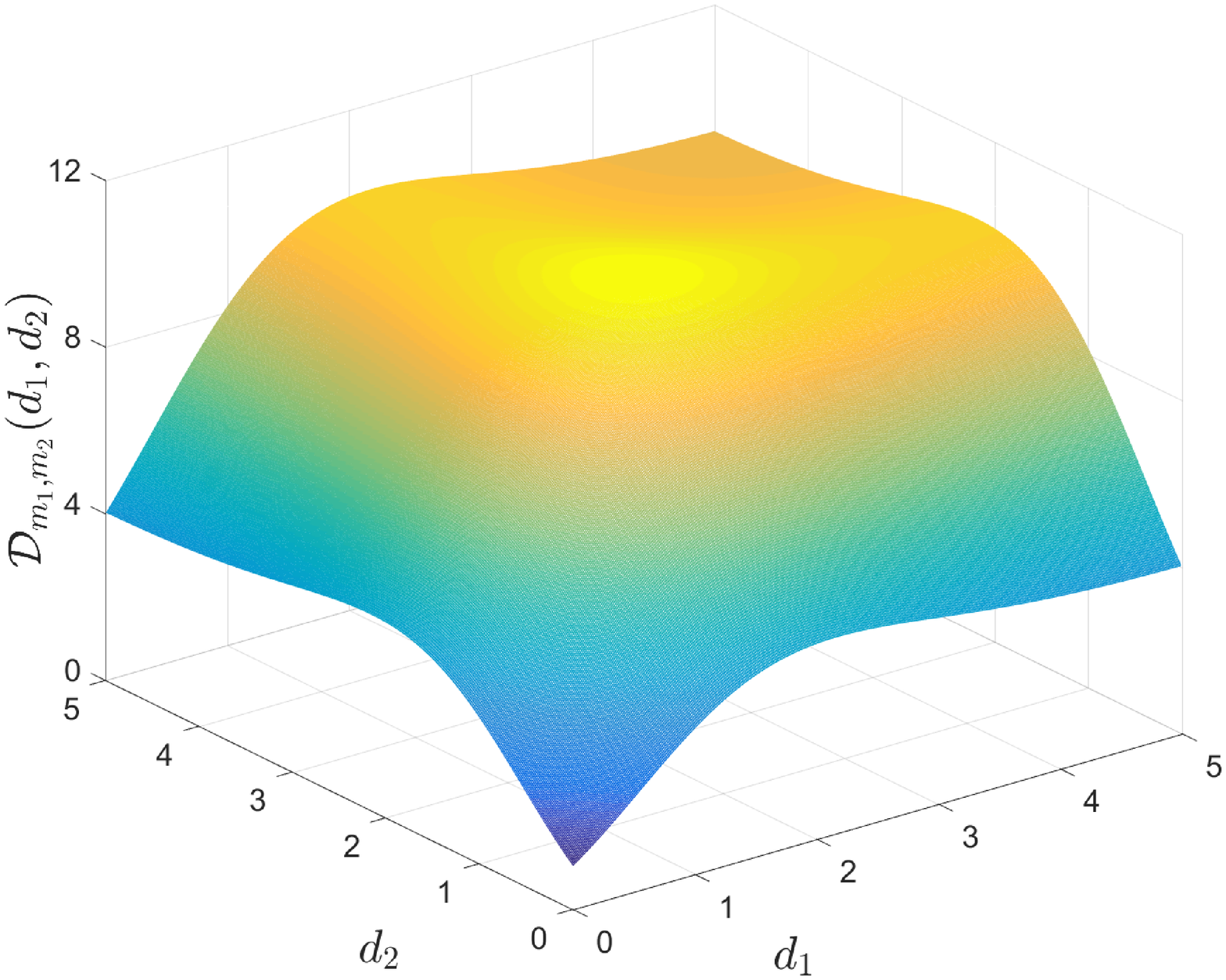,height=6cm, width=7.5 cm} \quad
\epsfig{file=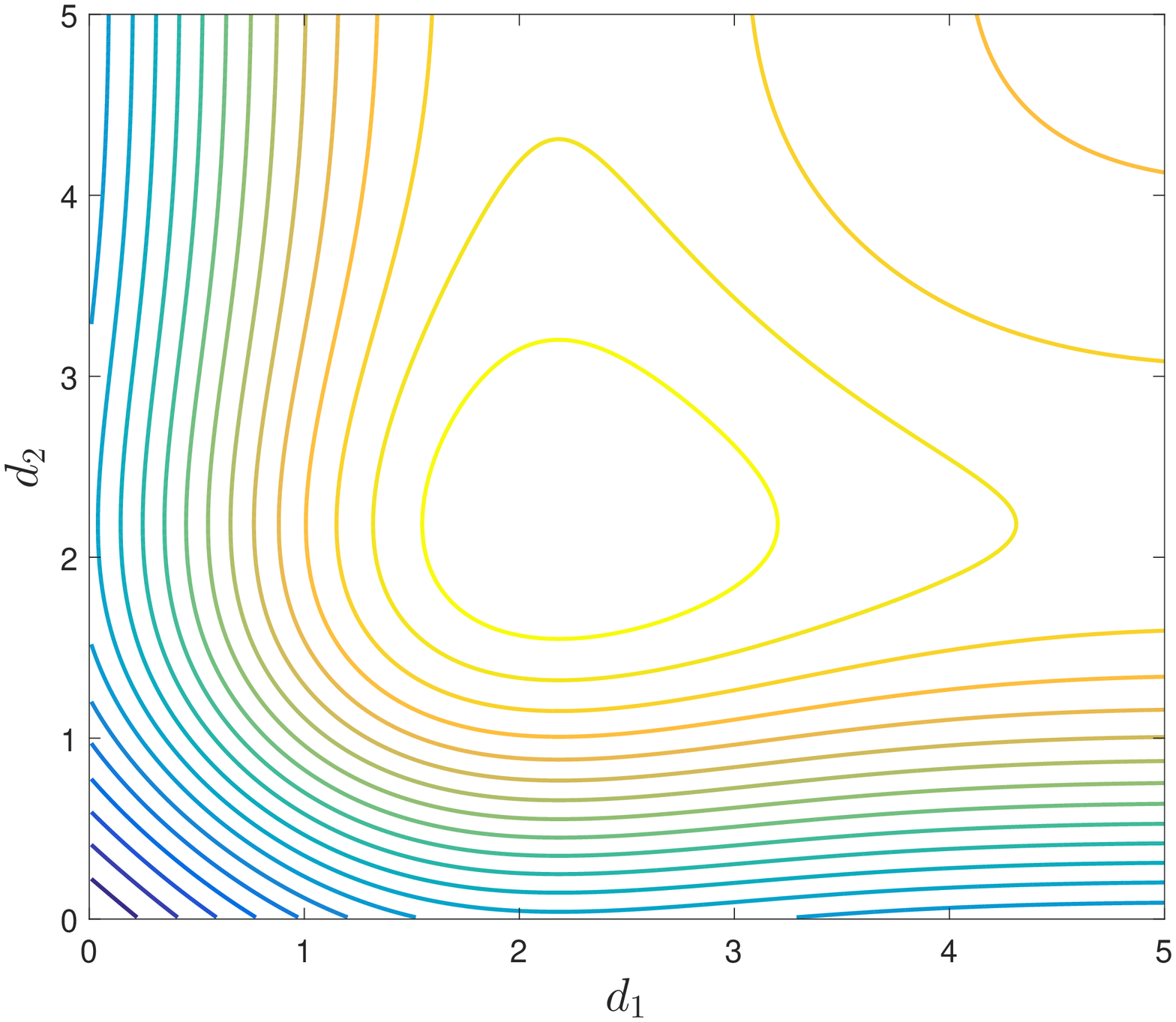,height=5.5 cm}}

\centerline{\hbox to 9 truecm {\scriptsize (a) \hfill (b)}}

\smallskip

\leavevmode
\hbox{
\epsfig{file=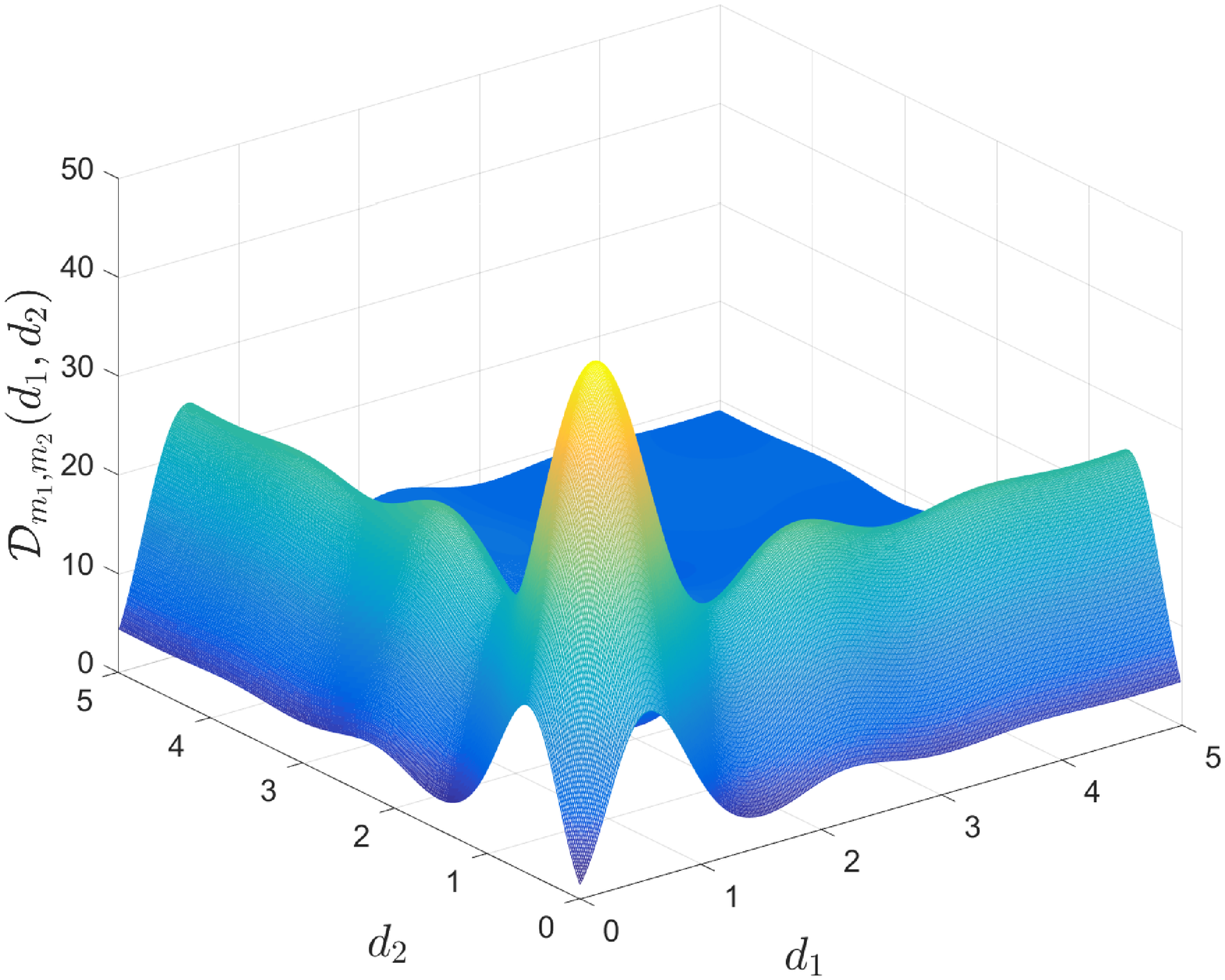,height=6cm, width=7.5 cm} \quad
\epsfig{file=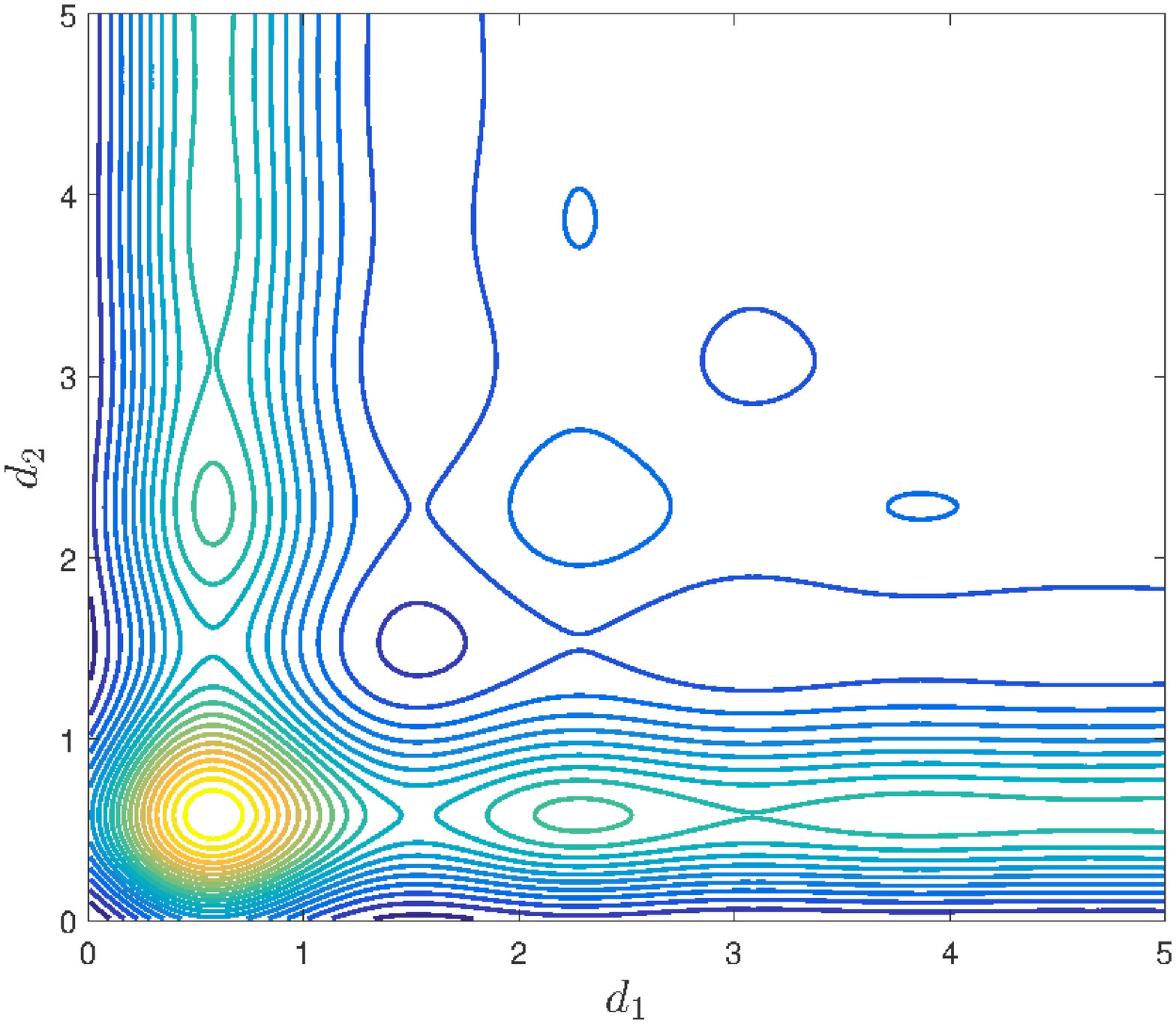,height=5.5cm}}

\centerline{\hbox to 9 truecm {\scriptsize (c) \hfill (d)}}

\end{center}
\caption{Bivariate objective function \ ${\mathcal D}_{m_1,m_2}$ \ for a three-point design for  (a) \ $\lambda=1, \ \omega=1$ \ and (c) \ $\lambda=1, \ \omega=4$, \ together with the corresponding contour plots (b) and (d), respectively. }
\label{fig1}
\end{figure}
Consider now the special case of a constant trend, that is the model
\begin{equation}
  \label{c_mod_s}
Z(t)=m+Y(t).
\end{equation}
In this situation \ $f_1(t)\equiv 1$ \ and \ $f_2(t)\equiv 0$, \ so 
\begin{equation*}
Z_1(t)=m_1+Y_1(t), \qquad Z_2(t)=m_2+Y_2(t),
\end{equation*} 
and the expression in \eqref{Qn_general} reduces to
\begin{equation}
   \label{Qn}
Q(n)\!=\!1\!+\!\sum_{\ell=1}^{n-1} g(d_{\ell}), \quad \text{where} \quad g(x)\!:=\!\frac{1\!-\!2{\mathrm e}^{-\lambda x}\cos(\omega x)\!+\!{\mathrm e}^{-2\lambda x}}{1\!-\!{\mathrm e}^{-2\lambda x}}, \ x\!>\!0, \ \ \text{and} \ \ g(0)\!:=\!0.
\end{equation}
Hence, in order to obtain the D-optimal design, one has to find
the maximum in \ $\boldsymbol d=(d_1,d_2,\ldots ,d_{n-1})$ \ of
\begin{equation}
  \label{Dopt_m}
{\mathcal D}_{m_1,m_2}(\boldsymbol d):=\det \big(\mathcal
I_{m_1,m_2}(n)\big)=\left (1+\sum_{\ell=1}^{n-1} g(d_{\ell})\right ) ^2.
\end{equation}

\begin{thm}
  \label{thm1}
Consider the complex model \eqref{c_mod_s} with \ $\omega\ne 0$ \ observed in design points  \ $0\leq t_1<t_2 < \ldots <t_n$. \ The D-optimal design for estimating the trend parameter is equidistant with \ $d_1=d_2=\ldots =d_{n-1}=d^*$, \ where \ $d^*$ \ is the (existing) global maximum point of \ $g(x)$.
\end{thm}

Observe that for \ $\omega=0$ \ we have 
\begin{equation*}
Q(n)=1+\sum_{\ell=1}^{n-1} \frac {1-{\mathrm e}^{-\lambda d_{\ell}}}{1+{\mathrm e}^{-\lambda d_{\ell}}},
\end{equation*}
which is exactly the Fisher information on the constant trend parameter of a shifted real valued stationary OU process with covariance parameter \ $\lambda$ \ \citep{ks,zba09}. In this case the D-optimal design on trend is also equidistant, however, with the increase of this equal distance the information is also increasing. According to the statement of Theorem \ref{thm1} this is not the case for the complex OU process as there exists an optimal distance which provides the highest information.

\begin{ex}
  \label{ex1}
As an illustration consider a three-point design. Figures \ref{fig1}a and \ref{fig1}c show  \ the bivariate objective function \ ${\mathcal D}_{m_1,m_2}$ \ for \ $\lambda=1, \ \omega=1$ \ and  \ $\lambda=1, \ \omega=4$ \ together with the corresponding contour plots (Figures \ref{fig1}b and \ref{fig1}d, respectively).

\end{ex}

\begin{figure}[t]
\begin{center}
\leavevmode
\epsfig{file=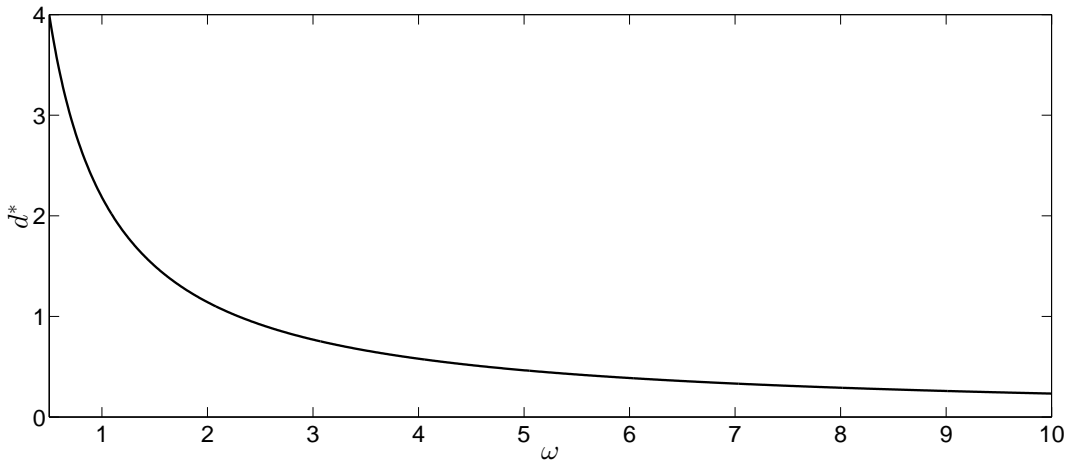,width=14 cm} 
\end{center}
\caption{Location \ $d^*$ \ of the global maximum of \ $g(x)$ \ as a function of the frequency \ $\omega$ \ for \ $\lambda =1$. }
\label{fig2}
\end{figure}

To get a better insight on the behavior of the optimal design, consider the function \ $g(x)$ \ defined by \eqref{Qn}. Figure \ref{fig2} shows the location \ $d^*$ \ of the global maximum of \ $g(x)$ \ as a function of the frequency \ $\omega$ \ for \ $\lambda =1$. \ The general case can obviously be obtained by rescaling both axes by the value of the damping parameter \ $\lambda$.

\subsection{Estimation of the covariance parameters}
  \label{subs:subs3.2}

Consider now the problem of estimating the damping parameter \ $\lambda$ \ and frequency \ $\omega$. \  Recalling again  the results of \citet{xia} and \citet{pazman}, the FIM on these parameters based on observations \ $\big\{ \big(Z_1(t_j),Z_2(t_j)\big), \ j=1,2,\ldots ,n\big\}$ \ has the form
\begin{equation}
  \label{FIM_cov}
\mathcal I_{\lambda,\omega}(n)=\begin{bmatrix}
        \mathcal I_{\lambda}(n) &  I_{\lambda,\omega}(n) \\
        I_{\lambda,\omega}(n) &  \mathcal I_{\omega}(n)
       \end{bmatrix},
\end{equation}
where
\begin{align}
\mathcal I_{\lambda}(n)&:=\frac 12 \tr \left\{C^{-1}(n)\frac{\partial C(n)}{\partial
    \lambda}C^{-1}(n)\frac{\partial C(n)}{\partial \lambda} \right\}, \nonumber \\
\mathcal I_{\omega}(n)&:=\frac 12 \tr \left\{C^{-1}(n)\frac{\partial C(n)}{\partial
    \omega }C^{-1}(n)\frac{\partial C(n)}{\partial \omega} \right\}, \label{FIM_cov_entry}\\
I_{\lambda,\omega}(n)&:=\frac 12 \tr \left\{C^{-1}(n)\frac{\partial
    C(n)}{\partial \lambda}C^{-1}(n)\frac{\partial C(n)}{\partial \omega} \right\}. \nonumber
\end{align}
Note, that here \ $\mathcal I_{\lambda}(n)$ \ and \ $\mathcal I_{\omega}(n)$ \ are Fisher information on parameters \ $\lambda$ \ and \ $\omega$, \ respectively, taking the other parameter as a nuisance.

\begin{thm}
  \label{thm2}
Consider the two-dimensional analogue of the complex model \eqref{c_mod_s} with \ $\omega\ne 0$ \ observed in design points  \ $0\leq t_1<t_2 < \ldots <t_n$. \ Then
\begin{equation}
   \label{FIM_cov_form}
\mathcal I_{\lambda}(n)=\sum_{\ell=1}^{n-1} \varphi(d_{\ell}), \qquad
\mathcal I_{\omega}(n)=\sum_{\ell=1}^{n-1} \psi(d_{\ell})
 \qquad \text{and} \qquad I_{\lambda,\omega}(n)=0, 
\end{equation}
where  \ $d_j=t_{j+1}-t_j, \ j=1,2,\ldots ,n-1$, 
\begin{equation}
   \label{PhiPsi}
\varphi(x)\!:=\!\frac{x^2{\mathrm e}^{-\lambda x}\cosh (\lambda x)}{\sinh^2(\lambda x)}, \quad  \psi(x)\!:=\!\frac{x^2{\mathrm e}^{-\lambda x}}{\sinh(\lambda x)}, \ \  x>0, \ \ \text{and} \ \ \varphi(0)\!:=\!\frac 1{\lambda ^2}, \ \ \psi(0)\!:=\!0.
\end{equation}
\end{thm}

Observe that none of the entries of the FIM \ ${\mathcal I}_{\lambda,\omega}(n)$ \ depends on the frequency parameter. Further, \  ${\mathcal I}_{\lambda}(n)/2$ \ coincides with the Fisher information on the covariance parameter of a real valued OU process given by \citet{zba09}. Note that here we consider two-dimensional OU processes, which justifies the halving of the information, however, due to this connection the first statement of the following theorem is a direct consequence of \citet[][Theorem 4.2]{zba09}.

\begin{thm}
   \label{thm3}
Consider the two-dimensional analogue of the complex model \eqref{c_mod_s} with \ $\omega\ne 0$ \ observed in design points  \ $0\leq t_1<t_2 < \ldots <t_n$. 
\begin{itemize}
  \item[a)] The D-optimal design for the damping parameter \ $\lambda$ \ maximizing the Fisher information \ ${\mathcal I}_{\lambda}(n)$ \ does not exist within the class of admissible designs.
 \item[b)] The D-optimal design for the frequency \ $\omega$ \ maximizing  the Fisher information \ ${\mathcal I}_{\omega}(n)$ \  is equidistant with
\begin{equation*}
d_1=d_2= \ldots =d_{n-1}=d_{\lambda}:=\frac 1{\lambda} \Big (W\big(-2{\mathrm e}^{-2}\big )/2+1\Big) \approx \frac {0.7968}{\lambda},
\end{equation*}
where \ $W(x)$ \ denotes the Lambert $W$-function \citep{cghjk}.
 \item[c)] The D-optimal design for both covariance parameters of the complex OU process maximizing \ $\det \big({\mathcal I}_{\lambda,\omega}(n)\big)$ \ is equidistant with
\begin{equation*}
d_1=d_2= \ldots =d_{n-1}=d^{\circ}/\lambda ,
\end{equation*}
where \ $d^{\circ}\approx 0.4930$ \ is the unique positive solution of
\begin{equation}
   \label{Dopt_eq}
1-d-2d{\mathrm e}^{-2d}-{\mathrm e}^{-4d}=0.
\end{equation}
\end{itemize}
\end{thm}

We remark that in the case of one-dimensional OU processes an admissible design for the covariance parameter \ $\lambda$ \ can be obtained by introducing the nugget effect, see e.g. \citet{shhs}.

\subsection{Estimation of all parameters}
  \label{subs:subs3.3}

\begin{figure}[t]
\begin{center}
\leavevmode
\hbox{
\epsfig{file=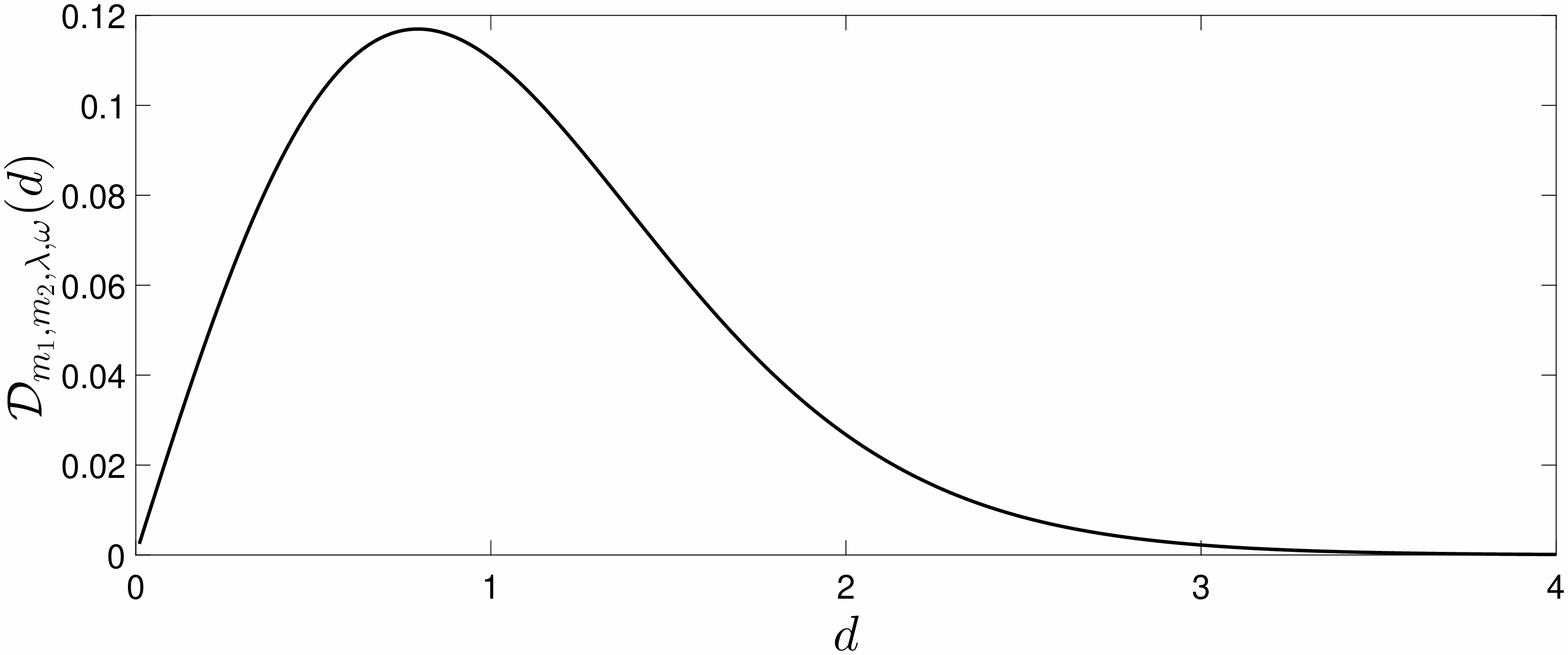,width=7.8cm} \quad
\epsfig{file=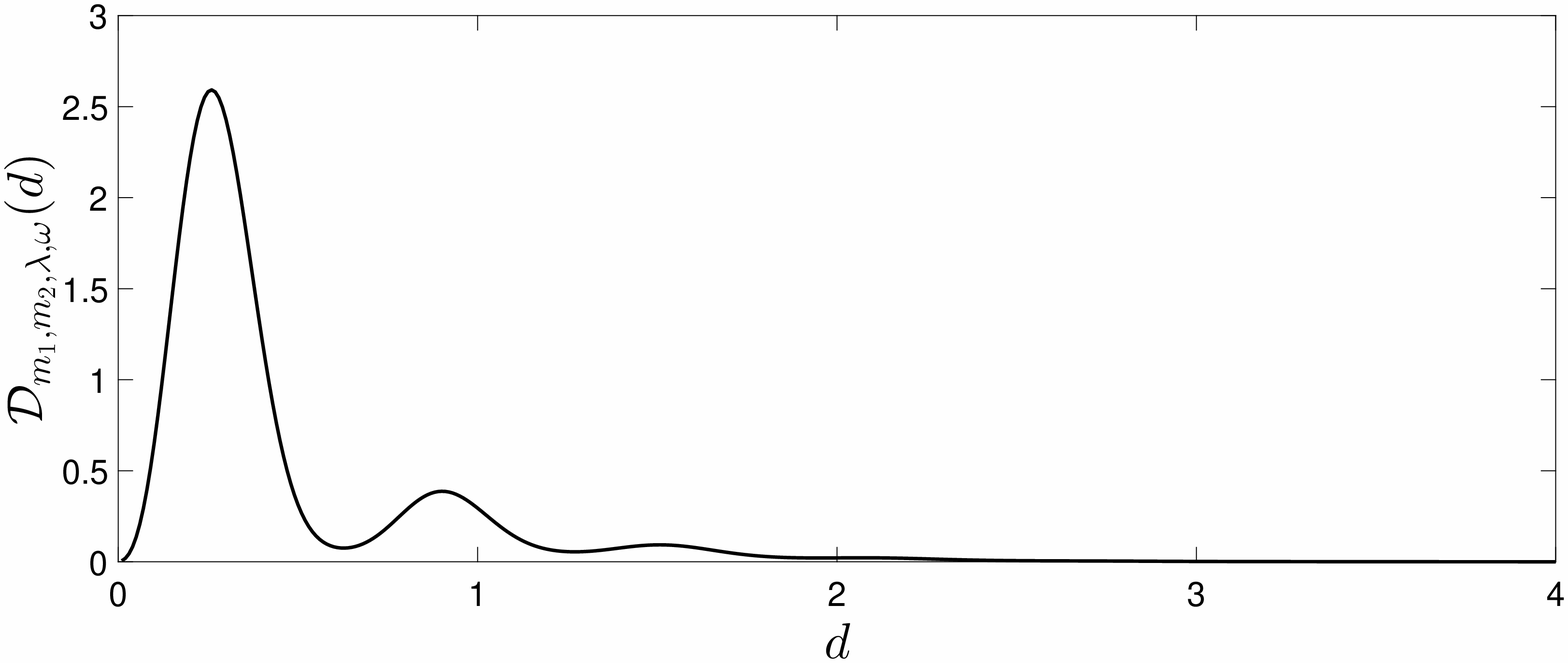,width=7.8cm}}

\centerline{\hbox to 9 truecm {\scriptsize (a) \hfill (b)}}
\end{center}
\caption{Objective function \ ${\mathcal D}_{m_1,m_2,\lambda,\omega}$ \ for a two-point design for  (a) \ $\lambda=1, \ \omega=1$ \ and (b) \ $\lambda=1, \ \omega=4$. }
\label{fig3}
\end{figure}
In the most general case one has to estimate both the components \ $m_1,m_2$ \ of the mean and the covariance parameters \ $\lambda$ \ and \ $\omega$. \ The FIM on these parameters equals
\begin{equation}
   \label{FIM_all}
\mathcal I_{m_1,m_2,\lambda,\omega}(n)=\begin{bmatrix}
        \mathcal I_{m_1,m_2}(n) &  {\mathbf 0}_{2,2} \\
         {\mathbf 0}_{2,2} &  \mathcal I_{\lambda,\omega}(n)
       \end{bmatrix},
\end{equation}
where \ $\mathcal I_{m_1,m_2}(n)$ \ and \ $\mathcal I_{\lambda,\omega}(n)$ \ are the information matrices defined in Lemma \ref{lem1} and by \eqref{FIM_cov}, respectively, and \ ${\mathbf 0}_{k,\ell}$ \ denotes the \ $k\times \ell$ \ matrix of zeros. Hence, according to the results of Lemma \ref{lem1} and Theorem \ref{thm2}, the D-optimal design for all four parameters maximizes in \ $\boldsymbol d=(d_1,d_2,\ldots ,d_{n-1})$ \ 
\begin{equation}
  \label{Dopt_all}
{\mathcal D}_{m_1,m_2,\lambda,\omega}(\boldsymbol d):=\det \big(\mathcal I_{m_1,m_2,\lambda,\omega}(n)\big)=\left (1+\sum_{\ell=1}^{n-1} g(d_{\ell})\right ) ^2\left (\sum_{\ell=1}^{n-1} \varphi (d_{\ell})\right )\left (\sum_{\ell=1}^{n-1} \psi(d_{\ell})\right ),
\end{equation}
where the functions \ $g$, \ $\varphi$ \ and \ $\psi$ \ are defined by \eqref{Qn} and \eqref{PhiPsi}.

\begin{figure}[t]
\begin{center}
\leavevmode
\hbox{
\epsfig{file=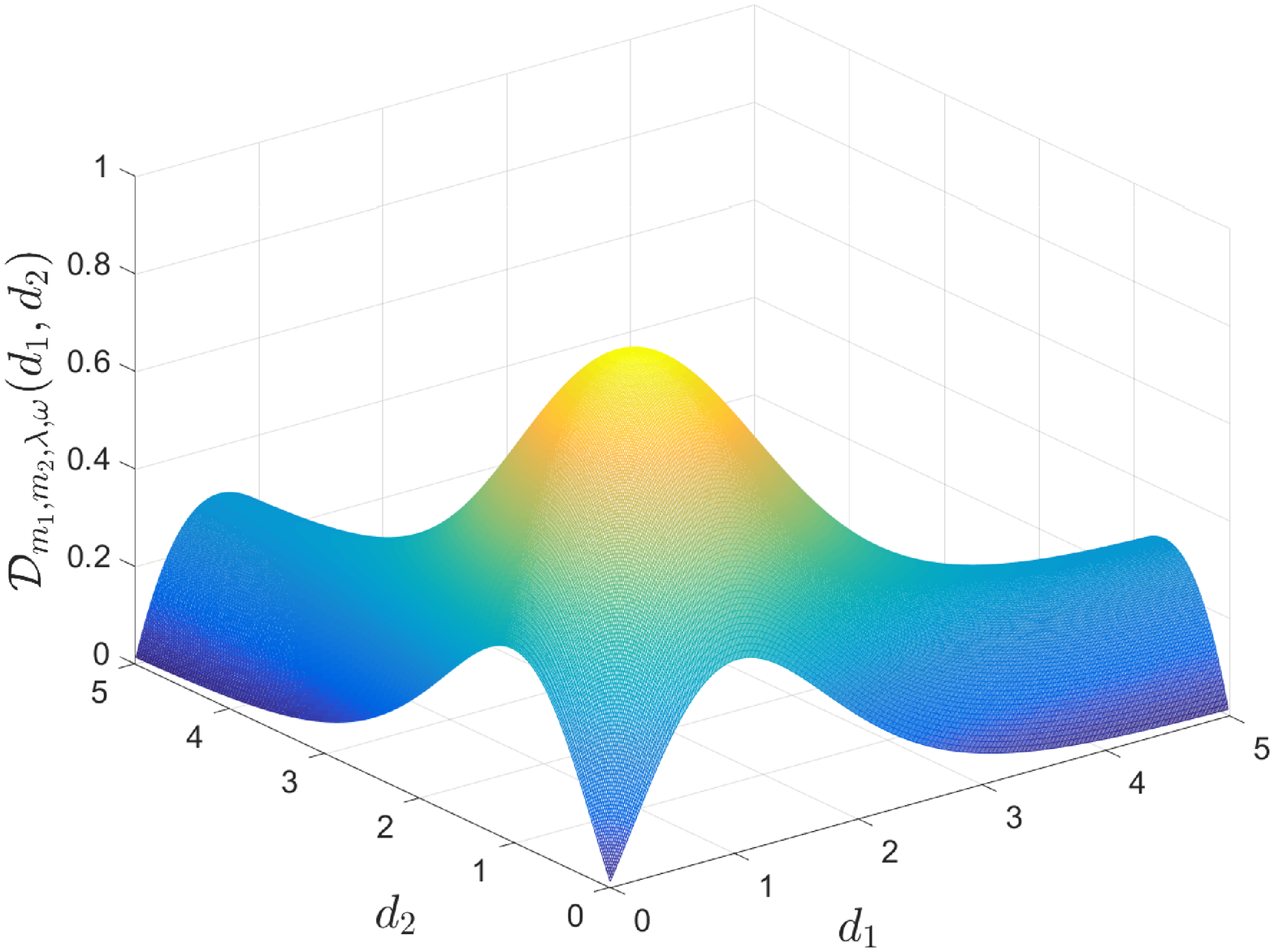,height=6cm, width=7.5 cm} \quad
\epsfig{file=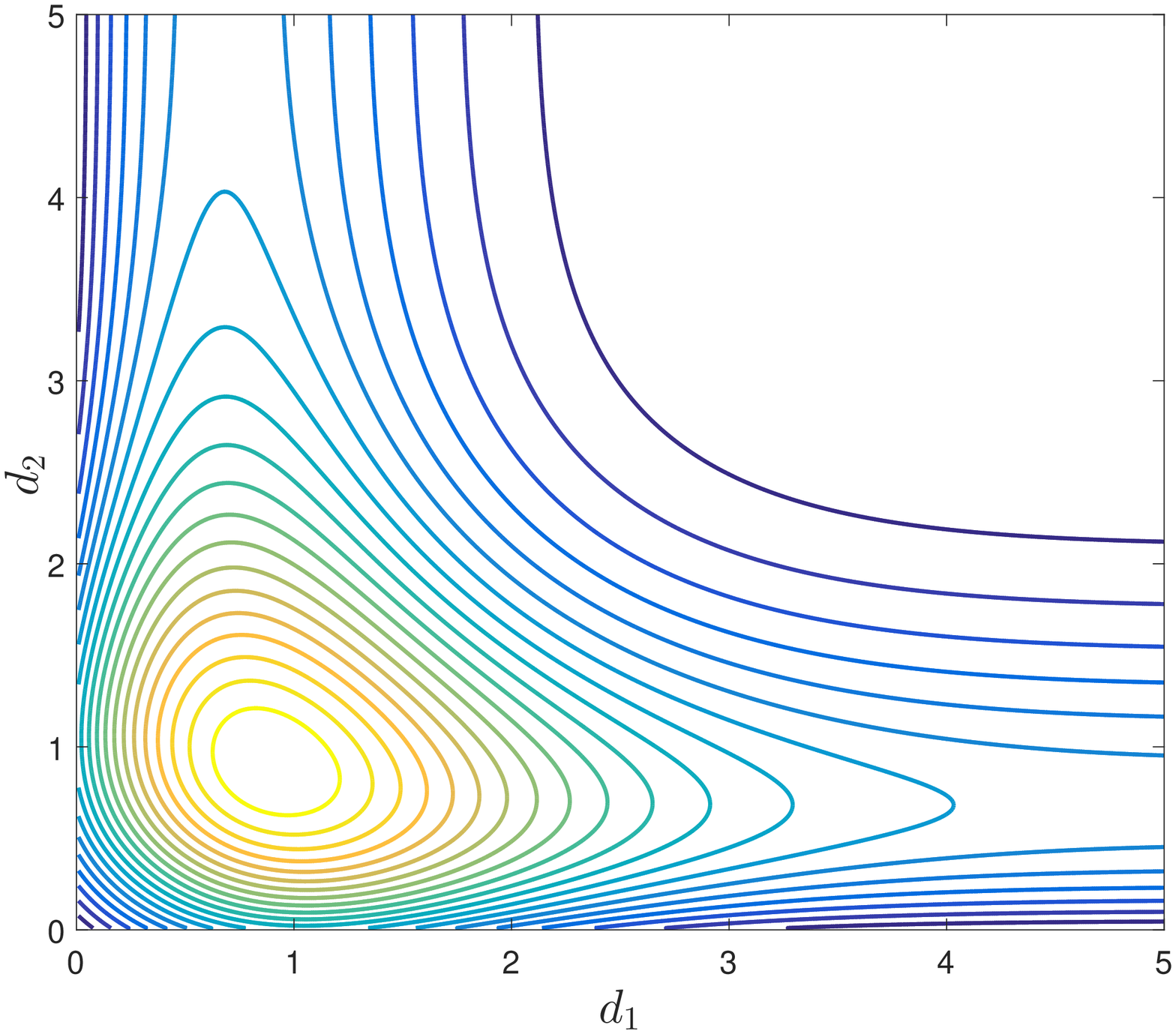,height=5.5 cm}}

\centerline{\hbox to 9 truecm {\scriptsize (a) \hfill (b)}}

\smallskip

\leavevmode
\hbox{
\epsfig{file=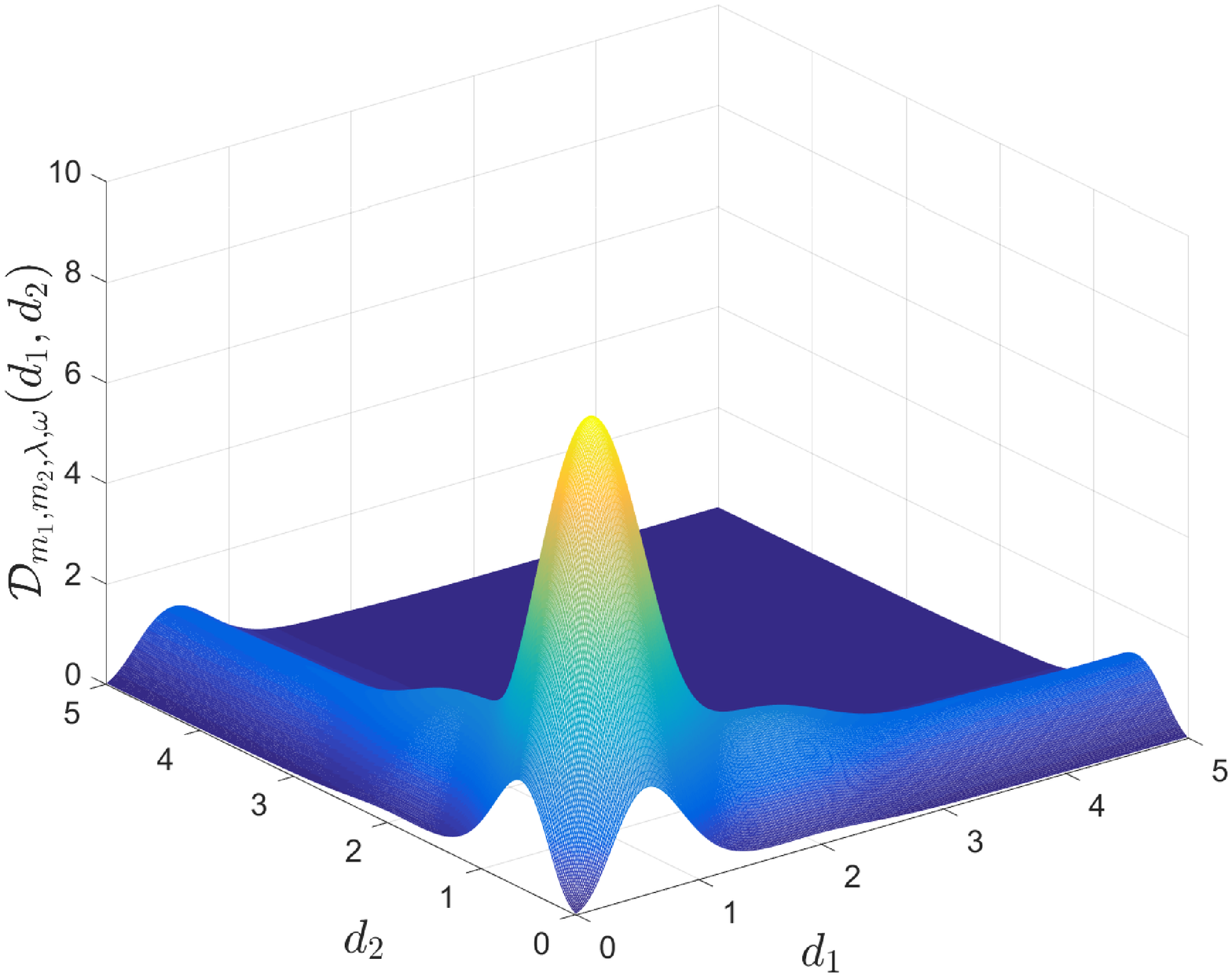,height=6cm, width=7.5 cm} \quad
\epsfig{file=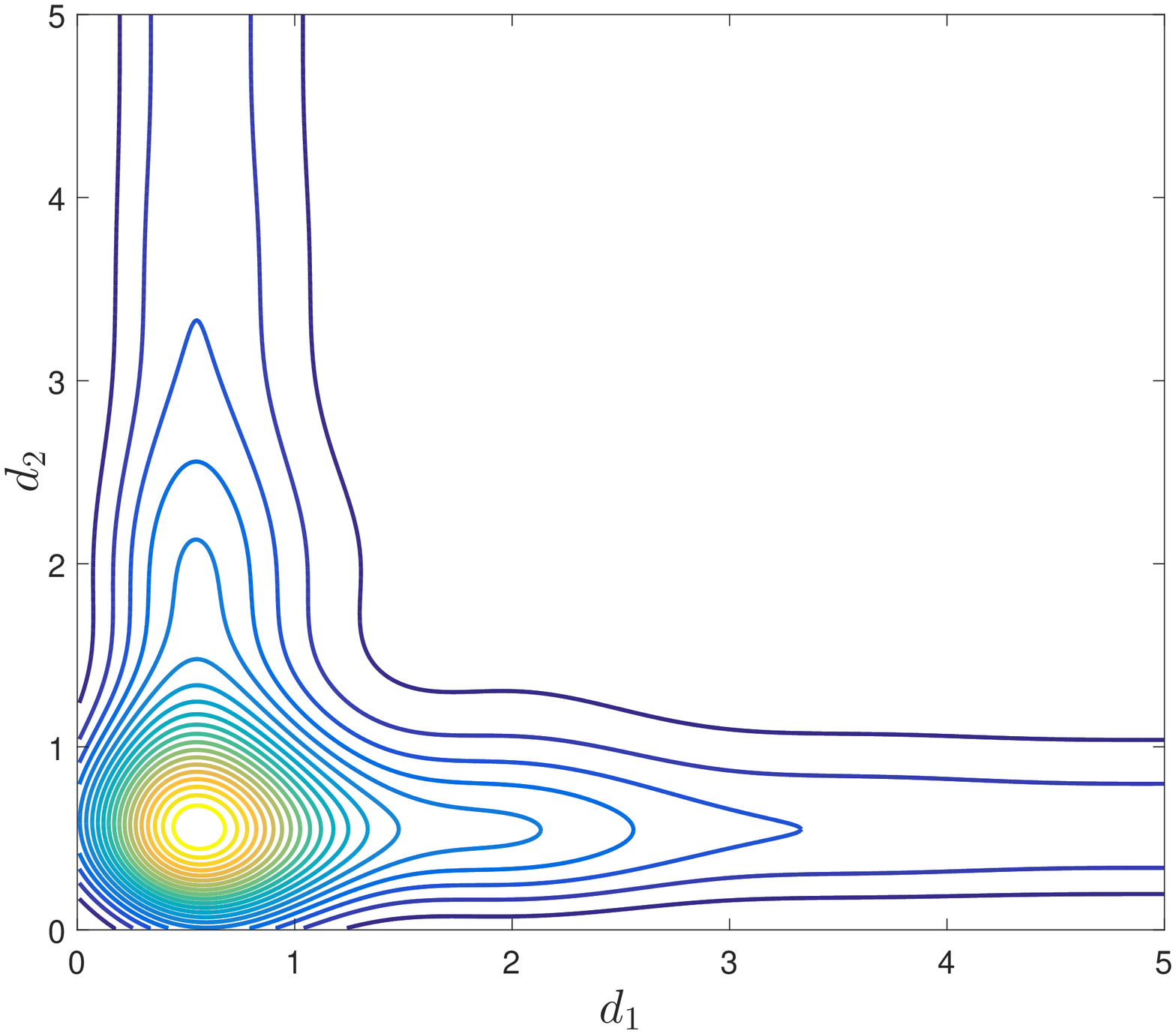,height=5.5cm}}

\centerline{\hbox to 9 truecm {\scriptsize (c) \hfill (d)}}

\end{center}
\caption{Bivariate objective function \ ${\mathcal D}_{m_1,m_2,\lambda,\omega}$ \ for a three-point design for  (a) \ $\lambda=1, \ \omega=1$ \ and (c) \ $\lambda=1, \ \omega=4$, \ together with the corresponding contour plots (b) and (d), respectively. }
\label{fig4}
\end{figure}

\begin{ex}
  \label{ex2}
Consider first the simplest case, that is a two-point design, when the objective function \eqref{Dopt_all} is univariate.  Figures \ref{fig3}a and \ref{fig3}b, showing the graph of \ ${\mathcal D}_{m_1,m_2,\lambda,\omega}(d)$ \ for \ $\lambda=1, \ \omega=1$ \ and  \ $\lambda=1, \ \omega=4$, \ respectively, clearly illustrate the existence of a finite optimal design.
\end{ex}

\begin{ex}
  \label{ex3}
Similar to Example \ref{ex1}, consider again a three-point design. In Figures \ref{fig4}a and \ref{fig4}c the bivariate objective function \ ${\mathcal D}_{m_1,m_2,\lambda,\omega}$ \ is plotted for \ $\lambda=1, \ \omega=1$ \ and  \ $\lambda=1, \ \omega=4$, \ together with the corresponding contour plots (Figures \ref{fig4}b and \ref{fig4}d, respectively).
\end{ex}  

\begin{figure}[t]
\begin{center}
\leavevmode
\epsfig{file=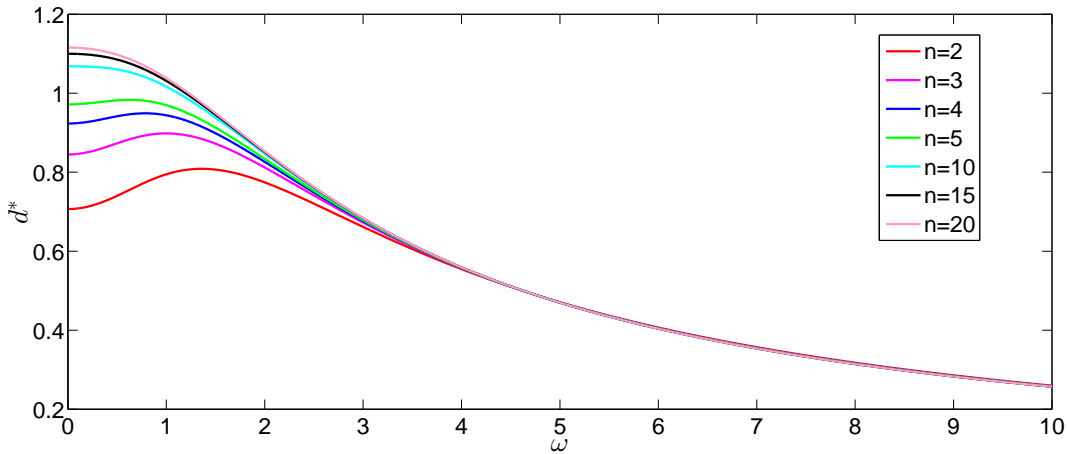,width=14 cm} 
\end{center}
\caption{Optimal distance \ $d^*$ \ of the equidistant $n$-point D-optimal design for estimation of all parameters as a function of the frequency \ $\omega$ \ for \ $\lambda =1$.}
\label{fig5}
\end{figure}

Although the objective function \ ${\mathcal D}_{m_1,m_2,\lambda,\omega}(\boldsymbol d)$ \ is too complicated for treating it analytically, numerical results in higher dimensions show that there exists a D-optimal design and it is equidistant. In Figure \ref{fig5} the optimal distance \ $d^*$ \ of the equidistant $n$-point D-optimal design is plotted as a  function of the frequency \ $\omega$ \ for \ $\lambda =1$. \ Again, the general case can easily be obtained by rescaling both axes by the damping \ $\lambda$. \ Note that the larger the frequency, the smaller the effect of the number of design points on the optimal design.

\section{Application}
  \label{sec:sec4}

The application of the obtained designs, especially the ones given in Theorem  \ref{thm3} can be applied in the assessment of the quality of parameter estimation for damping parameter \ $\lambda$ \ and frequency parameter \ $\omega$ \ in Kolmogorov's model \eqref{chandler} of the Chandler wobble. The 
maximum likelihood estimator and sufficient statistics for \  $\lambda$ \ are given in \citet{arato2,arato}. However, as noted after Theorem \ref{thm3}, for the drift parameter there exist no admissible design, unless we apply nugget effect.

For the model of \citet{aks} two different estimates of the damping parameter are given by \citet{pancenko} \ ($\widehat\lambda=0.3$) \ and \cite{wy}  \ ($\widehat\lambda=0.01$). \ By the second statement of Theorem \ref{thm3} the D-optimal design for frequency \ $\omega$ \ is equidistant with an optimal lag of \ $
d_{\lambda}\approx \frac{0.7968}{\lambda}$. \
However, the various estimates for \ $\lambda$ \ give a broad range of optimal equidistant times for measurements.

Further, due to the difference in the estimated values of the damping parameter \ $\lambda$, \ one might be interested in the sensitivity of the standardized setup \ ($\lambda=\omega=1$) \ with respect to the D-optimality criterion for estimation of the trend parameter. For the  $n$-point equidistant design this means the evaluation of \ $R(d,\lambda,1;n)$, \ where
\begin{equation}
 \label{Reff}
R(d,\lambda,\omega;n)\!:=\!\frac{\big(1\!+\!(n\!-\!1)g(d,1,1)\big)^2}{\big( 1\!+\!(n\!-\!1)g(d,\lambda ,\omega)\big)^2}, \quad \text{with} \quad g(x,\lambda,\omega)\!:=\!\frac{1\!-\!2{\mathrm e}^{-\lambda x}\cos(\omega x)\!+\!{\mathrm e}^{-2\lambda x}}{1\!-\!{\mathrm e}^{-2\lambda x}},
\end{equation} 
see also \eqref{Qn}. Here we analyze the situation for arbitrary \ $d$, \ since experimenter usually does not have a free choice of the lag, although for \ $\lambda=\omega=1$ \ the optimal value equals \ $d^*\approx 2.1835$ \ (see Theorem \ref{thm1}). Thus, this analysis incorporates all possible design spaces of form \ $[0,T_{\max}]$  \ where \ $T_{max}$ \ denotes the {\em upper bound\/} of the design space.

\begin{figure}[t]
\begin{center}
\leavevmode
\hbox{
\epsfig{file=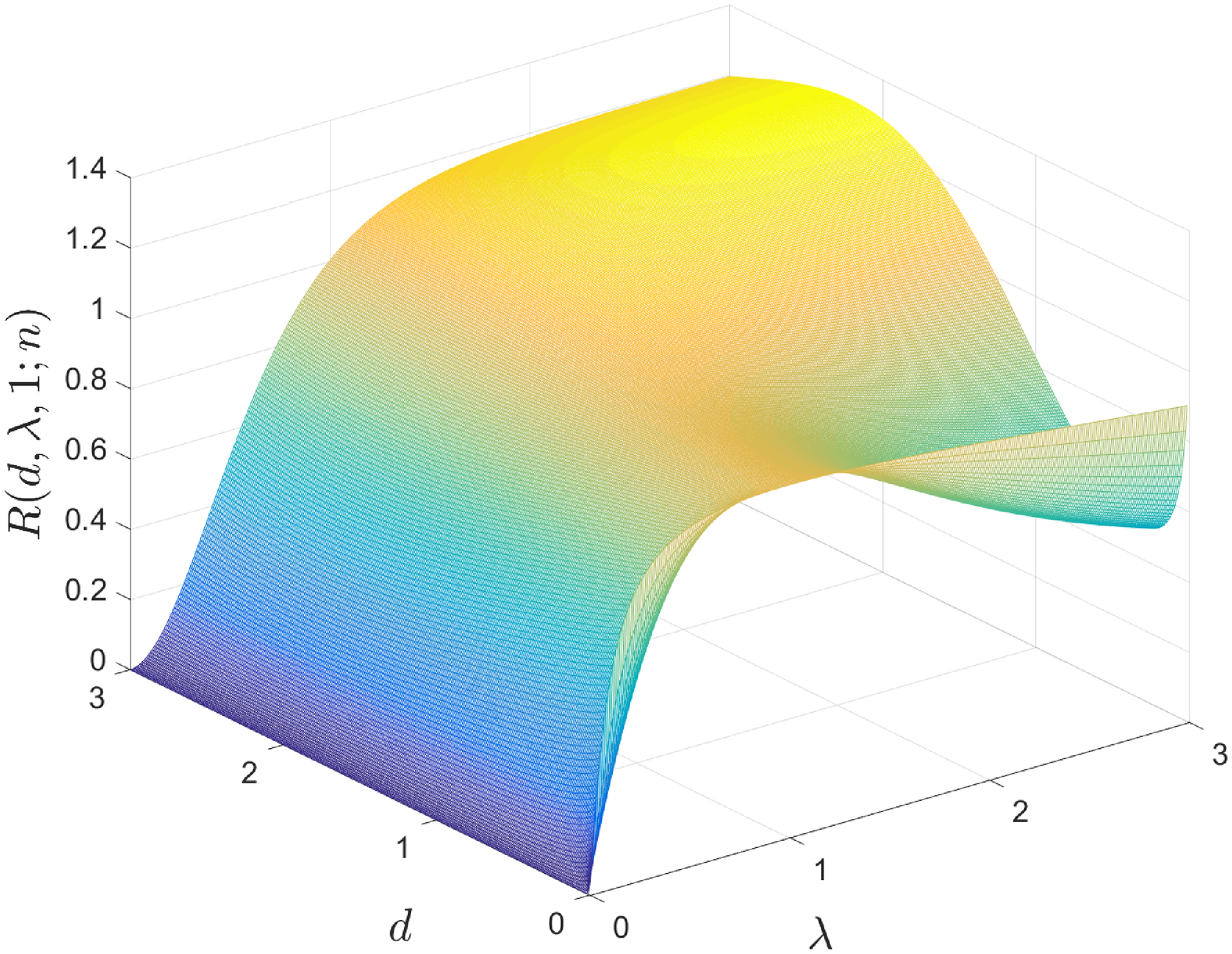,height=6cm, width=7.5 cm} \quad
\epsfig{file=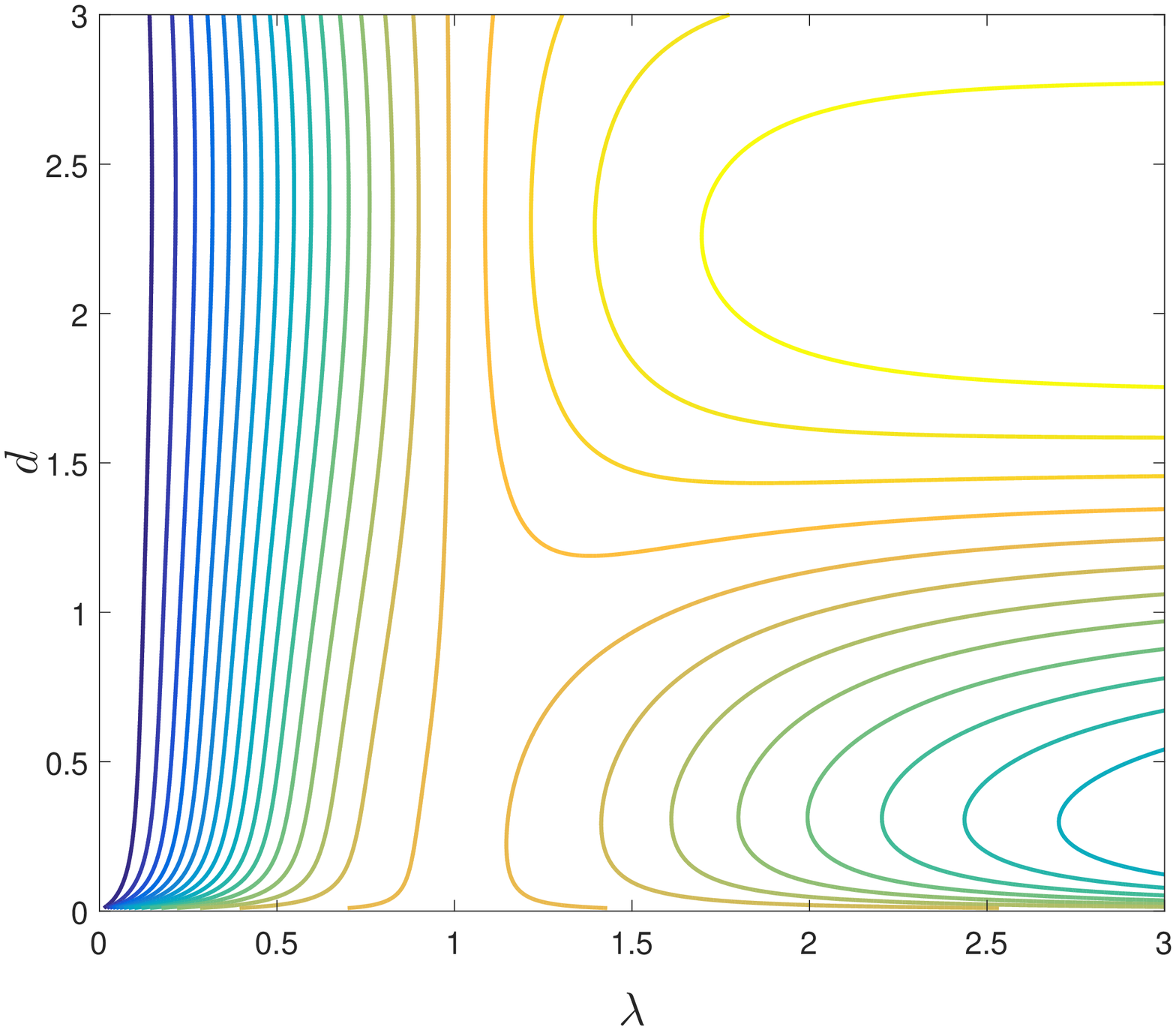,height=5.5 cm}}

\centerline{\hbox to 9 truecm {\scriptsize (a) \hfill (b)}}
\end{center}
\caption{Surface \ $R(d,\lambda,1;n)$  \ (a) and the corresponding contour plot (b) for a ten-point equidistant design.}
\label{fig6}
\end{figure}

The 3rd order Taylor expansion of \ $R(d,\lambda,1;n)$ \ around the origin \ $\lambda=0$ \ results in
\begin{equation*}
R(d,\lambda,1;n)=\frac{d^2\Big(1+(n-1)\frac{(1-2{\mathrm e}^{-d}\cos(d)+{\mathrm e}^{-2d})}{1-{\mathrm e}^{-2d}}\Big)^2}{(n-1)\big(\cos(d)-1\big)^2}\,\lambda^2+{\mathcal O}(\lambda^3),
\end{equation*}
which shows a very high sensitivity of the efficiency for the standard design 
with respect to small values of the damping parameter \ $\lambda$. \ The same phenomenon can be observed on Figure \ref{fig6} showing the surface \ $R(d,\lambda,1;n)$ \ and the corresponding contour plot for a ten-point equidistant design.

If one can use the optimal lag  \ $d^*\approx 2.1835$ \ for the standard design (i.e. has the possibility of choosing an arbitrary lag), then the 3rd order Taylor expansion of \ $R(d,\lambda,1;n)$ \ around the origin \ $\lambda=0$ \ reduces to
\begin{equation*}
R(d^*,\lambda,1;n)=\frac {1.9218\cdot \big(1+1.1569\cdot (n-1)\big)^2}{(n-1)^2}\lambda ^2+ {\mathcal O}(\lambda^3),
\end{equation*}
that is for \ $\lambda \approx 0$ \ we have
\begin{equation*}
  \lim_{n\to\infty}R(d^*,\lambda,1;n)\approx 2.5723 \lambda^2.
\end{equation*}
This also confirms the sensitivity of the efficiency for the standard design at small damping values \ $\lambda$.

\begin{figure}[t]
\begin{center}
\leavevmode
\hbox{
\epsfig{file=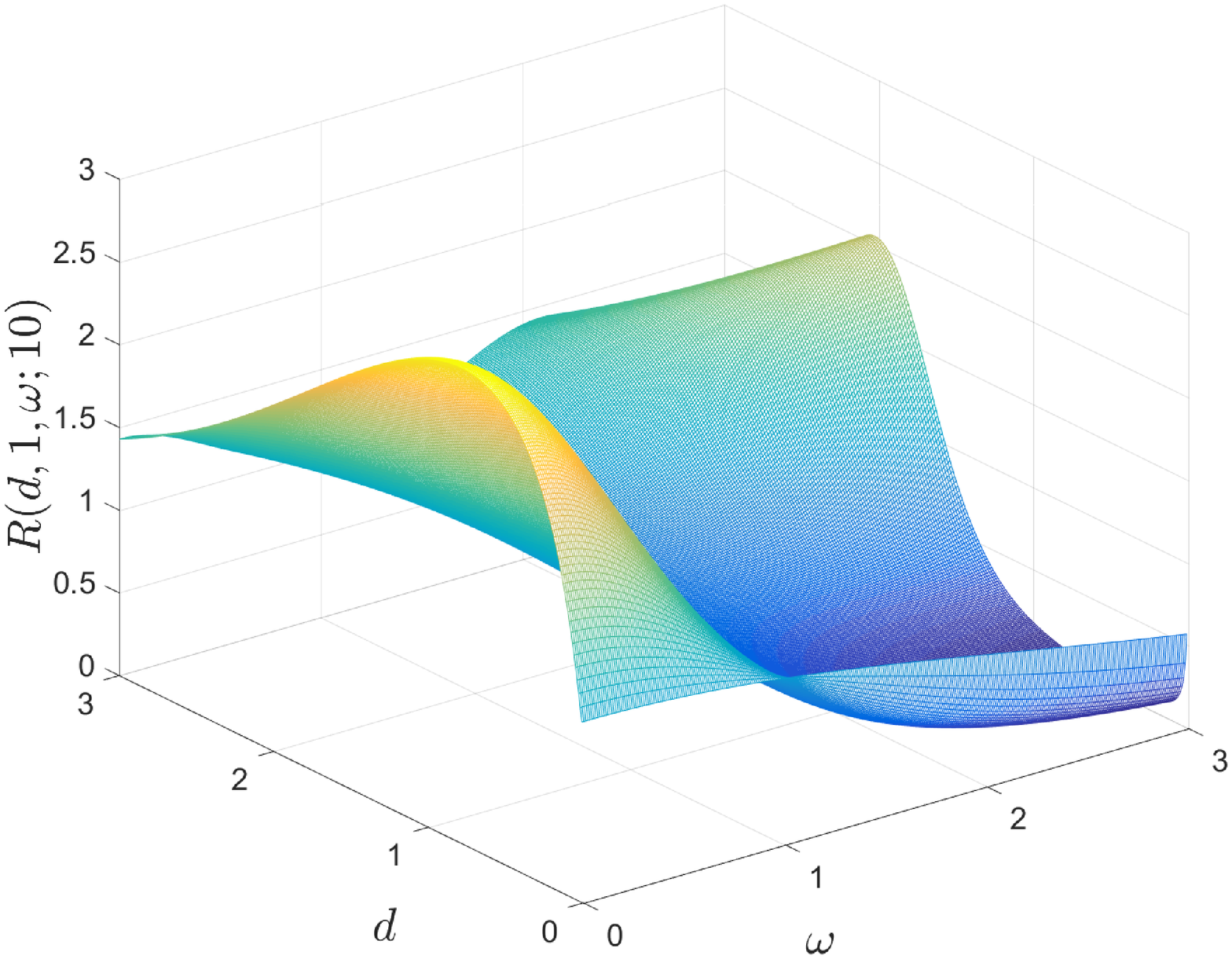,height=6cm, width=7.5 cm} \quad
\epsfig{file=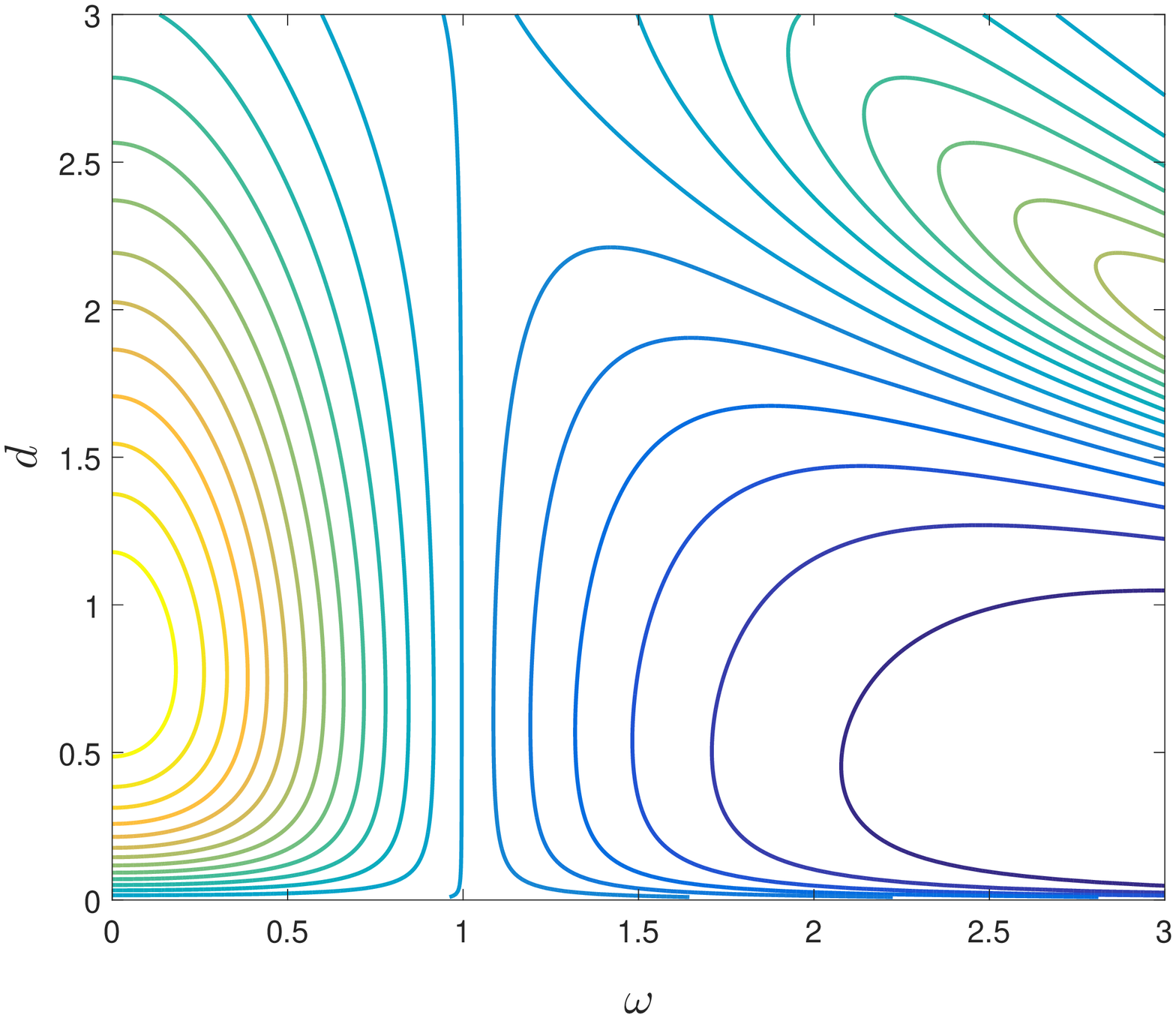,height=5.5 cm}}

\centerline{\hbox to 9 truecm {\scriptsize (a) \hfill (b)}}
\end{center}
\caption{Surface \ $R(d,1,\omega;n)$  \ (a) and the corresponding contour plot (b) for a ten-point equidistant design.}
\label{fig7}
\end{figure}

Consider now the sensitivity of efficiency for the standard design with respect to frequency \ $\omega$. \  For the  $n$-point equidistant design one has to evaluate \ $R(d,1,\omega;n)$ \ defined again by \eqref{Reff}, \ which for \ $\omega\approx 0$ \ indicates even higher sensitivity.  This can be derived from the form
\begin{align*}
R(d,1,\omega;n)=&\,\frac{\big((n-2){\mathrm e}^{-2d}-2(n-1){\mathrm e}^{-d}\cos(d)+n\big)^2}{\big((n-2){\mathrm e}^{-2d}-2(n-1){\mathrm e}^{-d}+n\big)^2}\\
&-\frac{2(n-1)d^2{\mathrm e}^{-d}\big((n-2){\mathrm e}^{-2d}-2(n-1){\mathrm e}^{-d}\cos(d)+n\big)^2}{\big((n-2){\mathrm e}^{-2d}-2(n-1){\mathrm e}^{-d}+n\big)^3}\,\omega^2+\mathcal O(\omega^4).
\end{align*}
of the 4th order Taylor expansion around \ $\omega=0$, \ and one can also observe it on Figure \ref{fig7}, where \ $R(d,1,\omega ;n)$ \ and the corresponding contour plot are shown for a ten-point equidistant design.

Again, for the optimal lag  \ $d^*\approx 2.1835$ \  the 4th order Taylor expansion of \ $R(d,1,\omega;n)$ \ around the origin \ $\omega=0$ \ reduces to
\begin{equation*}
R(d^*,1,\omega;n)=\frac{2.0455(9.0024n-1.2211)^2\big(6.2058n+1.5755-8.4655(n-1)\omega^2\big)}{(7.8777n+2)^3}+\mathcal O(\omega^4).
\end{equation*}
Hence, for \ $\omega\approx 0$ \ we have
\begin{equation*}
\lim _{n\to\infty}R(d^*,1,\omega;n)\approx 1.9476 - 2.3318\omega^2,
\end{equation*}
which also confirms the sensitivity of the efficiency at small frequencies \ $\omega$.

Recently \citet{mm} discovered that besides the well-known Chandler wobble phase jump in the 1920s, two other large phase jumps have been identified in the 1850s and 2000s. As in the 1920s, these phase jumps occurred contemporary with a sharp decrease in the Chandler wobble amplitude. 
However, sharp decrease of amplitude can drastically change the optimal design, as confirmed by Theorem \ref{thm3} for the case of a complex OU process.
This underpins the importance of further research on stochastic approach to Chandler wobble. Moreover, substantial relation of large seismic events to Chandler wobble excitation \citep{ocd} justify further studies of optimal designs for both damping and frequency parameters.

\subsection*{Acknowledgements}
This research was supported by the Hungarian -- Austrian intergovernmental S\&T cooperation program T\'ET\_{}15-1-2016-0046 and Project No. HU 11/2016. S\'andor Baran is grateful for the support of the J\'anos Bolyai Research Scholarship of the Hungarian Academy of Sciences. Milan Stehl\'\i k acknowledges the support of {\em FONDECYT Regular\/} Grant No 1151441.

\begin{appendix}
\section{Appendix}
  \label{sec:secA}

\subsection{Proof of Lemma \ref{lem1}}
   \label{subs:subsA.1}

Similar to the real valued case \citep{zba09}, the correlation matrix of observations \ $\big\{\! \big(Z_1(t_j),Z_2(t_j)\big), \ j=1,2,\ldots ,n\big\}$ \ of the two-dimensional process \eqref{2d_mod} equals
\begin{equation*}
C(n)=
     \begin{bmatrix}
        {\mathbb I}_2 & {\mathrm e}^{Ad_1} & {\mathrm e}^{A(d_1+d_2)} & 
        {\mathrm e}^{A(d_1+d_2+d_3)} &\dots &\dots & 
        {\mathrm e}^{A(\sum_{j=1}^{n-1}d_j)} \\
        {\mathrm e}^{A^{\top}\!d_1} &{\mathbb I}_2 &{\mathrm e}^{Ad_2} &
        {\mathrm e}^{A(d_2+d_3)} &\dots &\dots &{\mathrm e}^{A(\sum_{j=2}^{n-1}d_j)} \\
        {\mathrm e}^{A^{\top}\!(d_1+d_2)} &{\mathrm e}^{A^{\top}\!d_2} &{\mathbb I}_2 
        &{\mathrm e}^{Ad_3} &\dots &\dots &{\mathrm e}^{A(\sum_{j=3}^{n-1}d_j)} \\
         {\mathrm e}^{A^{\top}\!(d_1+d_2+d_3)} & {\mathrm e}^{A^{\top}\!(d_2+d_3)} 
         &{\mathrm e}^{A^{\top}\!d_3} &{\mathbb I}_2 &\dots &\dots &\vdots \\
        \vdots &\vdots &\vdots &\vdots &\ddots & &\vdots \\
        \vdots &\vdots &\vdots &\vdots & &\ddots &{\mathrm e}^{Ad_{n-1}} \\
        {\mathrm e}^{A^{\top}(\sum_{j=1}^{n-1}d_j)} &
        {\mathrm e}^{A^{\top}(\sum_{j=2}^{n-1}d_j)}
        &{\mathrm e}^{A^{\top}(\sum_{j=3}^{n-1}d_j)} &\dots &\dots &
        {\mathrm e}^{A^{\top}\!d_{n-1}} &{\mathbb I}_2
     \end{bmatrix}.
\end{equation*}
Short calculation shows that the inverse of \ $C(n)$ \ is given by
\begin{equation*}
  C^{-1}(n)=
     \begin{bmatrix}
        U_1 & -{\mathrm e}^{Ad_1}U_1 &0 & 0&\dots &\dots &0 \\
         -{\mathrm e}^{A^{\top}\!d_1}U_1 &V_2 &-{\mathrm e}^{Ad_2}U_2 &0 &\dots &\dots &0 \\
        0 & -{\mathrm e}^{A^{\top}\!d_2}U_2 &V_3  &-{\mathrm e}^{Ad_3}U_3&\dots &\dots &0 \\
        0 &0 &-{\mathrm e}^{A^{\top}\!d_3}U_3 &V_4&\dots &\dots &\vdots \\
        \vdots &\vdots &\vdots&\vdots  &\ddots & &\vdots \\
        \vdots &\vdots &\vdots &\vdots & &V_{n-1} &-{\mathrm e}^{Ad_{n-1}}U_{n-1} \\
        0 &0 &0 &\dots&\dots &-{\mathrm e}^{A^{\top}\!d_{n-1}}U_{n-1} &U_{n-1}
     \end{bmatrix},
\end{equation*}
where \ $U_k:=\big[{\mathbb I}_2-{\mathrm e}^{(A+A^{\top}\!)d_k}\big]^{-1},\ k=1,2,\dots,n-1$, \  and \ $V_k:=U_k+{\mathrm e}^{(A+A^{\top}\!)d_{k-1}}U_{k-1}, \ k=2,3,\dots,n-1$, \ which is a direct generalization of the result of \citet{ks} for the classical OU process.

Now, using \eqref{ou_cov} one can easily obtain 
\begin{alignat*}{2}
U_k&=\frac 1{1-{\mathrm e}^{-2\lambda d_k}}{\mathbb I}_2, \qquad {\mathrm e}^{Ad_k}U_k= \frac {{\mathrm e}^{-\lambda d_k}}{1-{\mathrm e}^{-2\lambda d_k}}\begin{bmatrix} \cos(\omega d_k)  & -\sin(\omega d_k) \\ \sin (\omega d_k) & \cos (\omega d_k) \end{bmatrix}, \quad &&k=1,2,\dots,n-1, \\
V_k&=\left (\frac 1{1-{\mathrm e}^{-2\lambda d_k}} +
\frac {{\mathrm e}^{-2\lambda d_{k-1}}}{1-{\mathrm e}^{-2\lambda d_{k-1}}}\right){\mathbb I}_2, \qquad &&k=2,3,\dots,n-1,
\end{alignat*}
which together with \eqref{ou_cov} specify both \ $C(n)$ \ and \ $C^{-1}(n)$. Finally, tedious but straightforward calculations lead us to \eqref{Qn_general}. \proofend

\subsection{Proof of Theorem \ref{thm1}}
   \label{subs:subsA.2}

Consider first the function \ $g(x)$ \ defined by \eqref{Qn}. As 
\begin{equation}
   \label{dg}
g^{\prime}(x)=\frac {r(x)}{\sinh^2(\lambda x)} \quad \text{with} \quad r(x):=\lambda \cosh(\lambda x)\cos(\omega x)\!+\!\omega \sinh (\lambda x)\sin(\omega x)\!-\!\lambda, \ x\!>\!0,
\end{equation}
the critical points of \ $g(x)$ \ are the roots of \ $r(x)$. \ Short calculation shows that \ $r(x)<0$ \ if and only if
\begin{equation}
  \label{dg_root}
\lambda \bigg[\tanh\Big(\frac{\lambda x}2\Big)\cot\Big(\frac{\omega x}2\Big)\bigg]^2 +2\omega \bigg[\tanh\Big(\frac{\lambda x}2\Big)\cot\Big(\frac{\omega x}2\Big)\bigg]-\lambda <0,
\end{equation}
and the quadratic polynomial in \eqref{dg_root} has two distinct roots 
\begin{equation*}
\varrho_-:=\frac {-\omega -\sqrt{\omega ^2+\lambda ^2}}\lambda <0 \qquad \text{and} \qquad \varrho_+:=\frac {-\omega +\sqrt{\omega ^2+\lambda ^2}}\lambda >0.
\end{equation*}
Now, if \ $\omega\ne 0$, \ by the properties of the cotangent and hyperbolic tangent functions there exist points \ $\big\{x_{n}^{-}, \ n\in{\mathbb N}\big\}$ \ such that \ $0<x_{n}^{-}<x_{n+1}^{-}$ \ and 
\begin{equation*}
\varrho_- <\tanh\Big(\frac{\lambda x_{n}^{-}}2\Big)\cot\Big(\frac{\omega x_{n}^{-}}2\Big)<\varrho_+, \qquad \text{that is} \qquad r\big(x_{n}^{-}\big)<0, \quad n\in{\mathbb N}.
\end{equation*}
Similar arguments prove the existence of points \ $\big\{x_{n}^{+}, \ n\in{\mathbb N}\big\}$ \ satisfying \ $0<x_{n}^{+}<x_{n+1}^{+}$ \ and \ $r\big(x_{n}^{+}\big)>0, \  n\in{\mathbb N}$. \ Thus, \ $r(x)$ \ (and \ $g^{\prime}(x)$) \ has an infinite number of sign changes, moreover, since \ $\lim_{x\to 0}g^{\prime}(x)=\frac {\omega ^2}{2\lambda }+\frac {\lambda }2>0$, \ the first change is from positive to negative. Denote by \ $\big\{x_{n}^{\circ}, \ n\in{\mathbb N}\big\}$ \ the roots of \ $r(x)$, that is the critical points of \ $g(x)$, \ 
where again, $0<x_{n}^{\circ}<x_{n+1}^{\circ}, \  n\in{\mathbb N}$, \ and \ $x_{1}^{\circ}$ \ is a local maximum of \ $g(x)$. \ Further, we have 
\begin{equation*}
g^{\prime\prime}\big(x_{n}^{\circ}\big)= \big(\lambda ^2 +\omega^2)\frac {\cos \big(\omega x_{n}^{\circ}\big)}{\sinh \big(\lambda x_{n}^{\circ}\big)}.
\end{equation*}
Assume that \ $g^{\prime\prime}\big(x_{n}^{\circ}\big)=0$, \ that is \ $\omega x_{n}^{\circ} =\frac {\pi}2+\pi k_n$ \ for some \ $k_n\in {\mathbb Z}$. \
In this case
\begin{equation}
  \label{ddg}
0=r\big(x_{n}^{\circ}\big)=\pm\omega \sinh \bigg(\frac {\lambda}{\omega} \Big(\frac {\pi}2+\pi k_n\Big) \bigg)-\lambda,
\end{equation}
where the positive sign will apply if \ $k_n$ \ is odd, and the negative sign if it is even. However, the Taylor series expansion of the hyperbolic sine implies
\begin{equation*}
\pm\omega \sinh \bigg(\frac {\lambda}{\omega} \Big(\frac {\pi}2+\pi k_n\Big) \bigg)
-\lambda =\pm\lambda \Big(\frac {\pi}2\mp 1+\pi k_n\Big)\pm\sum_{\ell=1}^{\infty} \frac {\lambda^{2\ell +1}}{(2\ell +1)!\,\omega^{2\ell}} \Big(\frac {\pi}2+\pi k_n\Big)^{2\ell+1}\ne 0,
\end{equation*}
which contradicts to \eqref{ddg}. Hence, for the critical points \ $\big\{x_{n}^{\circ}, \ n\in{\mathbb N}\big\}$ \ of \ $g(x)$ \ either  \ $g^{\prime\prime}\big(x_{n}^{\circ}\big)<0$ \ or  \ $g^{\prime\prime}\big(x_{n}^{\circ}\big)>0$, \ so they are either local maxima or local minima, respectively. This ensures the existence of a global maximum of \ $g$ \ as well at some point \ $d^*\in \big\{x_{n}^{\circ}, \ n\in{\mathbb N}\big\}$, \ where \ $g^{\prime}(d^*)=0$, \ and \ $g^{\prime\prime}(d^*)<0$.

Obviously, instead of \ ${\mathcal D}_{m_1,m_2}$ \ defined by \eqref{Dopt_m} it suffices to maximize 
\begin{equation*}
F(d_1,d_2,\ldots ,d_{n-1}):=Q(n)=1+\sum_{\ell=1}^{n-1}g\big(d_{\ell}\big)
\end{equation*}
Now, consider the $(n-1)$-dimensional vector \ ${\boldsymbol d}^*:=(d^*,d^*, \ldots ,d^*)^{\top }$. \ As 
\begin{equation*}
\frac {\partial F}{\partial d_k}(d^*,d^*, \ldots ,d^*)=g^{\prime}(d^*)=0, \qquad k=1,2,\ldots, n-1,
\end{equation*}
${\boldsymbol d}^*$ \ is a critical point of \ $F$. \ Further, as
\begin{equation*}
\frac {\partial^2 F}{\partial d_{\ell}\partial d_k}(d^*,d^*, \ldots ,d^*)=0 \qquad \text{if \ $k\ne\ell$,} \qquad  k,\ell=1,2,\ldots, n-1,
\end{equation*}
and
\begin{equation*}
\frac {\partial^2 F}{\partial d_k^2}(d^*,d^*, \ldots ,d^*)=g^{\prime\prime}(d^*)<0, \qquad k=1,2,\ldots, n-1,
\end{equation*}
the Hessian of  \ $F$ \ at \ ${\boldsymbol d}^*$ \ is negative definite. Hence, \ ${\boldsymbol d}^*$ \ is a maximum point of \ $F$ \ and since the components \ $g(d^*)$ \ are the largest possible, this maximum is global. \proofend

\subsection{Proof of Theorem \ref{thm2}}
   \label{subs:subsA.3}

Using the expressions of the covariance matrix \ $C(n)$ \ of observations and its inverse \ $C^{-1}(n)$ \ given in the proof of Lemma \ref{lem1} and \eqref{FIM_cov}, the formulae of \eqref{FIM_cov_form} can be verified by induction, similar to the proof of Theorem 2 of \citet{bs15}. As an example consider the Fisher information on damping parameter \ $\lambda$, \ where one has to show
\begin{equation}
  \label{InfLambda}
\mathcal I_{\lambda}(n):=\frac 12 \tr \left\{C^{-1}(n)\frac{\partial C(n)}{\partial
    \lambda}C^{-1}(n)\frac{\partial C(n)}{\partial \lambda} \right\}= \sum_{\ell=1}^{n-1}\frac{2d_{\ell}^2q_{\ell}^2(1+q_{\ell}^2)}{(1-q_{\ell}^2)^2},
\end{equation}
with \ $q_k:={\mathrm e}^{-\lambda d_k}, \ k=1,2, \ldots ,n-1$.

For \ $n=2$ \ equation \eqref{InfLambda} holds trivially. Assume also that
\eqref{InfLambda} is true for some \ $n$, \ and we are going to show it for
\ $n+1$. \  Let 
\begin{equation*}
\Delta (n)\!:=\!\big(-(d_1+d_2+\ldots+d_n)q_1q_2\ldots q_n B_{1,\ldots,n} ,
-(d_2+d_3\ldots+d_n)q_2q_3\ldots q_n B_{2,\ldots ,n}, \ldots , -d_nq_n B_{n}\big)^{\top}\!\!\!,
\end{equation*}
where
\begin{equation*}
B_{k,\ldots ,n}:=\begin{bmatrix}
        \cos\big(\omega(d_k+\ldots+d_n)\big)& -\sin\big(w(d_k+\ldots+d_n)\big)\\
         \sin\big(\omega(d_k+\ldots+d_n)\big) &  \cos\big(\omega(d_k+\ldots+d_n)\big)
       \end{bmatrix}, \qquad k=1,2,\ldots, n.
\end{equation*}
Using the representations of \ $C(n)$ \ and \ $C^{-1}(n)$ \ given in Section \ref{subs:subsA.1}, one can easily see
\begin{equation*}
\frac{\partial C(n\!+\!1)}{\partial \lambda}\!=\!\left[
  \begin{BMAT}{c.c}{c.c}
     \frac{\partial C(n)}{\partial \lambda} & \Delta(n)\\
     \Delta^{\top} (n) & {\mathbf 0}_{2,2}
  \end{BMAT}
\right] \ \ \text{and} \ \
 C^{-1}(n\!+\!1)\!=\!\left[
  \begin{BMAT}{c.c}{c.c}
      C^{-1}(n) & \Lambda_{1,2}(n)  \\
     \Lambda_{1,2}^{\top}(n)  & (1-q_n^2)^{-1} {\mathbb I}_2 
  \end{BMAT}
\right] \!+\!\left[
      \begin{BMAT}{c.c}{c.c}
      \Lambda_{1,1}(n) & {\mathbf 0}_{2n,2}\\
      {\mathbf 0}_{2,2n} & {\mathbf 0}_{2,2}
  \end{BMAT}
     \right]\!,
\end{equation*}
where \ ${\mathbf 0}_{k,\ell}$ \ denotes the \ $k\times \ell$ \ matrix of zeros and 
\begin{equation*}
\Lambda_{1,1}(n):=\left[
      \begin{BMAT}{c.c}{c.c}
      {\mathbf 0}_{2n-2,2n-2} & {\mathbf 0}_{2n-2,2} \\
      {\mathbf 0}_{2,2n-2} & \frac {q_n^2}{1-q_n^2} {\mathbb I}_2 
      \end{BMAT}
     \right] \qquad \text{and} \qquad
     \Lambda_{1,2}(n):=\left[
      \begin{BMAT}{c}{c.c}
      {\mathbf 0}_{2n-2,2} \\
      -\frac {q_n}{1-q_n^2} B_{n}
      \end{BMAT}
      \right].
\end{equation*}
In this way
\begin{equation*}
C^{-1}(n+1)\frac{\partial C(n+1)}{\partial \lambda}=\left[
  \begin{BMAT}{c.c}{c.c}
     C^{-1}(n)\frac{\partial C(n)}{\partial \lambda} + {\mathcal K}_{1,1}(n)&
     C^{-1}(n) \Delta(n)\\
      {\mathcal K}_{1,2}(n) & \frac{d_n q_n^2}{1-q_n^2}{\mathbb I}_2 
  \end{BMAT}
\right]+\left[
      \begin{BMAT}{c.c}{c.c}
      {\mathcal K}_{2,1}(n) & {\mathcal K}_{2,2}(n) \\
       {\mathbf 0}_{2,2n}&  {\mathbf 0}_{2,2}
  \end{BMAT}
     \right],
\end{equation*}
with
\begin{alignat*}{2}
{\mathcal K}_{1,1}(n)&:=\left[
      \begin{BMAT}{c}{c.c}
      {\mathbf 0}_{2n-2,2n} \\
      -\frac {q_n}{1\!-\!q_n^2}B_n\Delta^{\top}(n)
      \end{BMAT}
     \right], \qquad
{\mathcal K}_{1,2}(n):=&& \frac{1}{1\!-\!q_n^2}\Delta^{\top}(n)-\frac{q_n}{1\!-\!q_n^2}B_n^{\top} \left[\begin{BMAT}{c.c}{c}
		\Delta^{\top}(n-1)& {\mathbf 0}_{2,2}
		\end{BMAT}\right], \\
{\mathcal K}_{2,1}(n)&:=\left[
      \begin{BMAT}{c.c}{c.c}
      {\mathbf 0}_{2n-2,2n-2} & {\mathbf 0}_{2n-2,2} \\
      \frac {q_n^2}{1\!-\!q_n^2}\Delta^{\top}(n-1) & {\mathbf 0}_{2,2}
      \end{BMAT}
     \right], \qquad
&&{\mathcal K}_{2,2}(n):=\left[
      \begin{BMAT}{c}{c.c}
      {\mathbf 0}_{2n-2,2} \\
      -\frac {d_nq_n^3}{1\!-\!q_n^2}B_n
      \end{BMAT}
      \right].
\end{alignat*}
Hence,
\begin{align}
\mathcal I_{\lambda}(n\!+\!1)=&\,\mathcal I_{\lambda}(n)+\tr
\left\{C^{-1}(n)\frac{\partial C(n)}{\partial \lambda}{\mathcal
    K}_{1,1}(n) \right\}+\tr \big\{ C^{-1}(n) \Delta(n) {\mathcal K}_{1,2}(n) \big \}+ \frac 12 \tr \big\{ {\mathcal K}_{1,1}^2(n)
\big \} \nonumber \\ &
+ \frac{d_n^2 q_n^4}{(1\!-\!q_n^2)^2}+\tr\left\{C^{-1}(n)\frac{\partial C(n)}{\partial \alpha}{\mathcal K}_{2,1}(n)\right\} + \tr\big\{{\mathcal K}_{1,2}(n){\mathcal K}_{2,2}(n)\big\}.  \label{IndStep}
\end{align}
After long but straightforward calculations one can get
\begin{align*}
\tr \left\{C^{-1}(n)\frac{\partial C(n)}{\partial \lambda}{\mathcal K}_{1,1}(n) \right\}=&\,-\sum_{\ell=1}^{n-1}\frac{2d_{\ell}^2q_{\ell}^2q_{\ell+1}^2\cdot\ldots\cdot q_n^2}{(1-q_{\ell}^2)(1-q_n^2)}=-\tr\left\{C^{-1}(n)\frac{\partial C(n)}{\partial \lambda}{\mathcal K}_{2,1}(n)\right\},\\ 
\tr \left\{C^{-1}(n) \Delta(n) {\mathcal
    K}_{1,2}(n) \right\}=&\,\frac {2d_n^2q_n^2}{\left(1-q_n^2\right)^2}-\frac {2d_n^2q_n^4}{\left(1-q_n^2\right)^2}, \\
\frac12\tr \left\{{\mathcal K}_{1,1}^2(n)
\right\} =&\,\frac {d_n^2q_n^4}{(1-q_n^2)^2},  \qquad \qquad \qquad \ \
\tr\left\{{\mathcal K}_{1,2}(n){\mathcal K}_{2,2}(n)\right\}=\frac {2d_n^2q_n^4}{\left(1-q_n^2\right)^2},
\end{align*}
so \eqref{IndStep} implies 
\begin{equation*}
\mathcal I_{\lambda}(n+1)=\mathcal I_{\lambda}(n)+\frac{2d_n^2q_n^2(1+q_n^2)}{(1-q_n^2)^2},
\end{equation*}
which completes the proof of \eqref{InfLambda}. The other two statements of 
Theorem \ref{thm2} can be verified in a similar way.
\proofend

\subsection{Proof of Theorem \ref{thm3}}
   \label{subs:subsA.4}

As it has been mentioned, statement a) directly follows from Theorem 4.2 of \citet{zba09}.

Now, consider the function \ $\psi(x)$ \ defined by  \eqref{FIM_cov_form}. As 
\begin{equation*}
\psi^{\prime}(x)=\frac {4x{\mathrm e}^{-2\lambda x}\big(1-\lambda x -{\mathrm e}^{-2\lambda x}\big)}{\big(1-{\mathrm e}^{-2\lambda x}\big)^2},
\end{equation*}
$\psi(x)$ \ has a single critical point at \ $d^*:=\frac 1{\lambda} \big (W(-2{\mathrm e}^{-2})/2+1\big)$, \ moreover, \ $\psi^{\prime\prime}(d^*)<0$, \ so \ $d^*$ \ is a global maximum point. Hence, statement b) can be verified using the same arguments as in the proof of Theorem \ref{thm1}.

According to \eqref{FIM_cov_form}, the D-optimal design on \ $\lambda$ \ and \ $\omega$ \ maximizes
\begin{align*}
G(d_1,d_2,\ldots ,d_{n-1}):=&\left (\sum_{\ell=1}^{n-1}\varphi (d_{\ell})\right)\left (\sum_{\ell=1}^{n-1}\psi (d_{\ell})\right) \\
&=\left(\sum_{\ell=1}^{n-1}\frac{2d_{\ell}^{\,2} {\mathrm e}^{-2\lambda d_{\ell}}\big(1+{\mathrm e}^{-2\lambda d_{\ell}}\big)}{\big(1-{\mathrm e}^{-2\lambda d_{\ell}}\big)^2}\right)\left(\sum_{\ell=1}^{n-1}\frac{2d_{\ell}^{\,2} {\mathrm e}^{-2\lambda d_{\ell}}}{\big(1-{\mathrm e}^{-2\lambda d_{\ell}}\big)}\right), \nonumber
\end{align*}
and obviously, it suffices to consider the case \ $\lambda=1$. \ The critical points of \ $G$ \ are solutions of the system
\begin{equation}
  \label{diffG}
\frac{\partial G}{\partial d_k}(d_1,d_2,\ldots ,d_{n-1})=\varphi^{\prime}(d_k) \left (\sum_{\ell=1}^{n-1}\psi (d_{\ell})\right)+\psi^{\prime}(d_k) \left (\sum_{\ell=1}^{n-1}\varphi (d_{\ell})\right)=0, 
\end{equation}
$k=1,2,\ldots ,n-1$, \ where
\begin{equation*}
\varphi^{\prime}(x)=\frac {4x{\mathrm e}^{-2\lambda x}\big(1-\lambda x -3\lambda x{\mathrm e}^{-2\lambda x}-{\mathrm e}^{-4\lambda x}\big)}{\big(1-{\mathrm e}^{-2\lambda x}\big)^3} <0, \qquad x>0.
\end{equation*}
Hence, for \ $\lambda=1$ \ the system of equations \eqref{diffG} is equivalent to
\begin{equation}
  \label{Gcritical}
\kappa(d_k)=\left (\sum_{\ell=1}^{n-1}\psi (d_{\ell})\right) \bigg/ \left (\sum_{\ell=1}^{n-1}\varphi (d_{\ell})\right), \qquad  k=1,2,\ldots ,n-1,
\end{equation}
where
\begin{equation*}
\kappa(x):=\frac{1-x +(x-2){\mathrm e}^{-2x}+{\mathrm e}^{-4x}}{x-1 +3x{\mathrm  e}^{-2x}+{\mathrm e}^{-4x}}, \qquad x>0.
\end{equation*}
Since \ $\kappa(x)$ \ is strictly monotone decreasing and continuous with \ $\lim_{x\to 0}\kappa(x)=1$ \ and $\lim_{x\to \infty}\kappa(x)=-1$, \ the solution of \eqref{Gcritical} should satisfy \ $d_1=d_2=\ldots =d_{n-1}=:d$. \ Hence, the critical points of \ $G$ \ have equal coordinates and at a critical point \eqref{Gcritical} reduces to \ $\kappa(d)=\big(1-{\mathrm e}^{-2d}\big)/\big(1+{\mathrm e}^{-2d}\big)$, \ which is equivalent to \eqref{Dopt_eq}. The latter equation has a unique positive solution \ $d^{\circ}$, so  \ $d_1=d_2=\ldots =d_{n-1}=d^{\circ}$ \ is the only critical point of \ $G$. \ Short calculation shows that at this point the Hessian of \ $G$ \ equals
\begin{align*}
H&=(n-1)\big[\varphi^{\prime\prime}(d^{\circ})\psi(d^{\circ})+\psi^{\prime\prime}(d^{\circ})\varphi(d^{\circ})\big]{\mathbb I}_{n-1}+2\big[\varphi^{\prime}(d^{\circ}) \psi^{\prime}(d^{\circ})\big]{\boldsymbol 1}_{n-1}{\boldsymbol 1}_{n-1}^{\top} \\
&\approx -0.5083\cdot (n-1){\mathbb I}_{n-1}-0.2754\cdot {\boldsymbol 1}_{n-1}{\boldsymbol 1}_{n-1}^{\top}, 
\end{align*}
where \ ${\boldsymbol 1}_k,\ k\in{\mathbb N},$ \ denotes the $k$-dimensional vector of ones. This means that for all \ ${\boldsymbol 0}\ne {\boldsymbol v}=(v_1,v_2,\ldots ,v_{n-1})^{\top}\in{\mathbb R}^{n-1}$ \ we have 
\begin{equation*}
{\boldsymbol v}^{\top}H{\boldsymbol v}=(n-1)\big[\varphi^{\prime\prime}(d^{\circ})\psi(d^{\circ})+\psi^{\prime\prime}(d^{\circ})\varphi(d^{\circ})\big] {\boldsymbol v}^{\top}{\boldsymbol v} + 2\big[\varphi^{\prime}(d^{\circ}) \psi^{\prime}(d^{\circ})\big]\left(\sum_{\ell=1}^{n-1}v_{\ell}\right)^2<0,
\end{equation*}
so the Hessian is negative definite. Hence, the unique critical point of \ $G$  \ is a global maximum, which completes the proof. \proofend

\end{appendix}

\end{document}